\documentclass{amsart}
\usepackage{amssymb}
\usepackage{amsfonts}
\usepackage{amsmath}
\usepackage{psfig}
\usepackage{psfrag}
\usepackage{graphicx}

\newtheorem{theorem}{Theorem}[section]
\newtheorem{lemma}[theorem]{Lemma}
\newtheorem{proposition}[theorem]{Proposition}
\newtheorem{definition}[theorem]{Definition}
\newtheorem{example}[theorem]{Example}

\newtheorem{remark}[theorem]{Remark}

\numberwithin{equation}{section}

\begin{document}

\title[Verification of mixing properties]{Verification of mixing properties in two-dimensional shifts of finite type}

\author{Jung-Chao Ban}
\address{Department of Applied Mathematics, National Dong Hwa University, Hualien 97401, Taiwan}
\email{jcban@mail.ndhu.edu.tw}
\thanks{The first author would like to thank the National Science Council, R.O.C. (Contract No. NSC 100-2115-M-259-009-MY2) and the National Center for Theoretical Sciences for partially supporting this research.}

\author{Wen-Guei Hu}
\address{Department of Applied Mathematics, National Chiao Tung University, Hsinchu 30010, Taiwan}
\email{wghu@mail.nctu.edu.tw}
\thanks{The second author would like to thank the National Science Council, R.O.C. and
the ST Yau Center for partially supporting this research.}

\author{Song-Sun Lin}
\address{Department of Applied Mathematics, National Chiao Tung University, Hsinchu 300, Taiwan}
\email{sslin@math.nctu.edu.tw}
\thanks{The third author would like to thank the National Science Council, R.O.C. (Contract No. NSC 103-2115-M-009-004) and
the ST Yau Center for partially supporting this research.}

\author{Yin-Heng Lin}
\address{Department of Applied Mathematics, National Chiao Tung University, Hsin-Chu 30010, Taiwan}
\email{amandalin@cathayholdings.com.tw}

\begin{abstract}
The degree of mixing is a fundamental property of a dynamical system. General multi-dimensional shifts cannot be systematically determined. This work introduces constructive and systematic methods for verifying the degree of mixing, from topological mixing to strong specification (or strong irreducibility) for two-dimensional shifts of finite type. First, transition matrices on infinite strips of width $n$ are introduced for all $n\geq 2$. To determine the primitivity of the transition matrices, connecting operators are introduced to reduce the order of high-order transition matrices to yield lower-order transition matrices. Two sufficient conditions for primitivity are provided; they are invariant diagonal cycles and primitive commutative cycles of connecting operators. After primitivity is established, the corner-extendability and crisscross-extendability are used to demonstrate topological mixing. In addition, the hole-filling condition yields the strong specification. All mentioned conditions can be verified to apply in a finite number of steps.
\end{abstract}

\maketitle

\section{Introduction}

\label{intro}
\hspace{0.5cm}

Multi-dimensional shift spaces represent an important and highly active area of research into topological dynamical systems. Such shifts also closely related to lattice models that are used in the scientific modeling of spatial structures. More precisely, when the lattice dynamical systems or coupled map lattices are spatial invariant and their equilibria assume only finite many values, the set of all stationary solutions forms a multi-dimensional shift space \cite{14,15,16,17}. Related investigations have been performed in statistical physics, chemistry \cite{4,6,7,10-1,15,18-1,24,25,26,26-1,27,36,37,38,39,40,47-1,48,48-1}, biology \cite{8,9}, image processing and pattern recognition \cite{14,16,17,21,22,23,30}.

Lattice models would be better understood if multi-dimensional shifts of finite type were better understood. The most interesting properties of shifts include their spatial entropy and various mixing properties, such as topological mixing and strong specification (or strong irreducibility). These properties are some of the important properties of dynamical systems \cite{11,12,13,13-1,13-2,22,23,27,29-1,35,40-1,40-2,42,42-1,43,45,47,50}. However, determining whether a given system exhibits topological mixing or strong specification in multi-dimensions is not easy. The intrinsic difficulty is related to the undecidability of the multi-dimensional coloring problem \cite{10,18,21,29,31,44,46,49}. For example, the extendability of local patterns on a finite lattice to a global pattern on $\mathbb{Z}^{2}$ is undecidable. Therefore, the mixing property of an arbitrary multi-dimensional shift cannot be determined. Nevertheless, this work provides some easily checked sufficient conditions for topological mixing and strong specification of certain two-dimensional shifts of finite type, which satisfy some non-degeneracy conditions of transition matrices $\mathbb{H}_{2}$ and $\mathbb{V}_{2}$. See Definition \ref{definition:3.3} and Theorem \ref{theorem:3.5}.

Let $\mathbb{Z}^{2}$ be a two-dimensional planar lattice. Vertex (or corner) coloring is considered first. For any $m,n\geq 1$ and $(i,j)\in\mathbb{Z}^{2}$, the $m\times n$ rectangular lattice with the left-bottom vertex $(i,j)$ is denoted by

\begin{equation*}
\mathbb{Z}_{m\times n}((i,j))=\left\{(i+n_{1},j+n_{2})\mid 0\leq n_{1}\leq m-1,0\leq n_{2}\leq n-1 \right\}.
\end{equation*}
In particular,
\begin{equation*}
\mathbb{Z}_{m\times n}=\mathbb{Z}_{m\times n}((0,0)).
\end{equation*}
Let $\mathcal{S}_{p}$ be a set of $p$ ($\geq 2$) colors (or symbols). For $m,n\geq 1$, $\Sigma_{m\times n}(p)=\mathcal{S}_{p}^{\mathbb{Z}_{m\times n}}$ is the set of all $m\times n$ local patterns or rectangular blocks.

Let $\mathcal{B}\subset\Sigma_{2\times 2}(p)$ be a basic set of admissible local patterns. For any lattice $R\subset\mathbb{Z}^{^{2}}$, the set of all $\mathcal{B}$-admissible patterns on $R$ is defined as

\begin{equation*}
\Sigma_{R}(\mathcal{B})=\left\{U\in\mathcal{S}_{p}^{R}: U\mid _{\mathbb{Z}_{2\times 2}((i,j))}\in\mathcal{B} \text{ if } \mathbb{Z}_{2\times 2}((i,j))\subset R
\right\}.
\end{equation*}
Let $\Sigma_{m\times n}(\mathcal{B})=\Sigma_{\mathbb{Z}_{m\times n}}(\mathcal{B})$ for $m,n\geq 2$.
$\Sigma(\mathcal{B})=\Sigma_{\mathbb{Z}^{2}}(\mathcal{B})$ is the set of all global patterns that can be constructed from the admissible local patterns in $\mathcal{B}$.

Traditionally, the admissible local patterns are specified on sublattices $\mathbb{Z}_{2\times 1}$ and $\mathbb{Z}_{1\times 2}$ with symbols in $\mathcal{S}_{p}$. Our approach to the two-dimensional problem begins with a study of infinite strips $\mathbb{Z}_{\infty\times n}$ and $\mathbb{Z}_{m\times\infty}$. Based on $\Sigma_{2\times n}(\mathcal{B})$ and $\Sigma_{m\times 2}(\mathcal{B})$, the transition matrices $\mathbb{H}_{n}$ on $\mathbb{Z}_{2\times n}$ and $\mathbb{V}_{m}$ on $\mathbb{Z}_{m\times 2}$ are introduced, and these apply to admissible patterns on $\mathbb{Z}_{\infty\times n}$ and $\mathbb{Z}_{m\times\infty}$, respectively.
Carefully arranging the local patterns on $\mathbb{Z}_{2\times 2}$ into the ordering matrices $\mathbf{X}_{2}$ and $\mathbf{Y}_{2}$, yields recursive formulae for $\mathbb{H}_{n}$ in $n$ and $\mathbb{V}_{m}$ in $m$, which are crucial in computing the spatial entropy \cite{1,2} and studying the mixing problem herein, as elucidated in Section 2. Notably, any $\mathbb{Z}^{2}$-shift of finite type can be represented by some $\mathcal{B}\subset\Sigma_{2\times 2}(p)$ for some $p\geq 2$ \cite{41}. Accordingly, only the case of $\mathcal{B}\subset\Sigma_{2\times 2}(p)$, $p\geq 2$, is considered here.

First, topological mixing is introduced. For any shift $\Sigma$ and any subset $R\subset \mathbb{Z}^{2}$, the restriction map is $\Pi_{R}(\Sigma):\Sigma\rightarrow \mathcal{S}_{p}^{R}$. Denote by $d$ the Euclidean metric on $\mathbb{Z}^{2}$. A $\mathbb{Z}^{2}$ shift $\Sigma$ is topologically mixing (or mixing, for short) if for any finite subsets $R_{1}$ and $R_{2}$ of $\mathbb{Z}^{2}$, a constant $M(R_{1},R_{2})$ exists such that for all $\mathbf{v}\in \mathbb{Z}^{2}$ with $d(R_{1},R_{2}+\mathbf{v})\geq M$, and for any two admissible patterns $U_{1}\in\Pi_{R_{1}}(\Sigma)$ and $U_{2}\in\Pi_{R_{2}+\mathbf{v}}(\Sigma)$, there exists a global pattern $W\in\Sigma$ with $\Pi_{R_{1}}(W)=U_{1}$ and $\Pi_{R_{2}+\mathbf{v}}(W)=U_{2}$; see \cite{50}.

$\Sigma$ has strong specification if a number $M(\Sigma)\geq 1$ exists such that for any two admissible patterns $U_{1}\in\Pi_{R_{1}}(\Sigma)$ and $U_{2}\in\Pi_{R_{2}}(\Sigma)$ with $d(R_{1},R_{2})\geq M$, where $R_{1},R_{2}$ are subsets of $\mathbb{Z}^{2}$, there exists a global pattern $W\in\Sigma$ with $\Pi_{R_{1}}(W)=U_{1}$ and $\Pi_{R_{2}}(W)=U_{2}$ \cite{50}. Clearly, strong specification implies topological mixing.

Some known results verify that $\Sigma(\mathcal{B})$ is topologically mixing or has strong specification \cite{3,42,42-1,48-2}.
Previously, in an investigation of pattern generation problems \cite{2}, the present authors introduced connecting operators to study the entropy of $\Sigma(\mathcal{B})$. In this work, connecting operators are also utilized to provide sufficient conditions for the topological mixing or strong specification of $\Sigma(\mathcal{B})$.

 A non-negative matrix $A$ is called primitive (or $N$-primitive) if there exists $N\geq 1$ such that each entry of $A^{n}$ is positive for all $n\geq N$. A matrix $A$ is called weakly primitive (or weakly $N$-primitive) if there exists $N\geq 1$ such that each entry of $A^{n}$ is positive except in positions of $A$ where a zero row or zero column is present for all $n\geq N$. The local crisscross-extendability and locally corner-extendable conditions are introduced in Section 3.

The main theorem for topological mixing is proven as Theorem \ref{theorem:3.14} and stated as follows.
\begin{theorem}
\label{theorem:1.1}
If
\begin{enumerate}
\item[(i)] $\mathcal{B}\subset \Sigma_{2\times 2}(p)$ is locally crisscross-extendable, and

\item[(ii)] $\mathcal{B}$ satisfies three of the locally corner-extendable conditions $C(i)$, $1\leq i\leq 4$,
\end{enumerate}
then
$\mathbb{H}_{n}(\mathcal{B})$ and $\mathbb{V}_{n}(\mathcal{B})$ are weakly primitive for all $n\geq 2$ if and only if
$\Sigma(\mathcal{B})$ is topologically mixing.
\end{theorem}

To provide checkable sufficient conditions for primitivity of $\mathbb{H}_{n}(\mathcal{B})$ and $\mathbb{V}_{n}(\mathcal{B})$, two sufficient conditions for the primitivity of $\mathbb{H}_{n}$ and $\mathbb{V}_{n}$ are introduced in Section 4; they are

\begin{enumerate}

\item[(i)] invariant diagonal cycles of connecting operator and

\item[(ii)] primitive commutative cycles of connecting operator.

\end{enumerate}
The invariant diagonal cycles look like a periodic structure of connecting operators, which is imposed to prove that for some given $q\geq 1$, the primitivity of $\mathbb{H}_{n+kq}(\mathcal{B})$ can be established from the primitivity of $\mathbb{H}_{n+(k-1)q}(\mathcal{B})$, $k\geq 1$; the conditions are then used to establish inductively  that $\mathbb{H}_{n}$ is primitive for all $n\geq 2$. When either condition (i) or condition (ii) applies, then only conditions concerning $\mathbb{H}_{n}$, $2\leq n\leq q+1$, have to be verified to ensure that $\mathbb{H}_{n}$ is primitive for all $n\geq 2$. More precisely, when $S$-invariant diagonal cycle $\overline{\beta}_{q}=\beta_{1}\beta_{2}\cdots\beta_{q}\beta_{1}$ of order $(m,q)$, with its invariant index set $\mathcal{K}$, exists, only the primitivity of $\underset{l\in\mathcal{K}}{\sum} H_{m,n;\beta_{1}}^{(l)}$ has to be verified for $2\leq n\leq q+1$. A similar result holds for primitive commutative cycles. See Section 4 for the details of the notation used; see Theorems 4.4 and 4.8 for detailed results.

Next, strong specification is considered. Strong specification is stronger than topological mixing. The hole-filling condition (HFC) introduced in Definition 5.1 is useful to provide checkable sufficient conditions for strong specification. HFC is closely related to the extension property called square filling \cite{40-1,40-2}. The main theorem for strong specification is Theorem \ref{theorem:6.4} and stated as follows.

Let $A=[a_{i,j}]_{n\times n}$ be a non-negative matrix; the index set of non-zero rows of $A$ and the index set of non-zero columns of $A$ are denoted by

\begin{equation}\label{eqn:1.51}
\begin{array}{ccc}
r(A)=\left\{i \mid \underset{j=1}{\overset{n}{\sum}}a_{i,j}>0\right\} & \text{and} &
c(A)=\left\{j \mid \underset{i=1}{\overset{n}{\sum}}a_{i,j}>0\right\},
\end{array}
\end{equation}
respectively.

\begin{theorem}
\label{theorem:1.2}
Given $\mathcal{B}\subset\Sigma_{2\times 2}(p)$, if there exists $k\geq 2$ such that
\begin{enumerate}
\item[(i)]  $r(\mathbb{H}_{k})=c(\mathbb{H}_{k})$ and $r(\mathbb{V}_{k})=c(\mathbb{V}_{k})$,

\item[(ii)] $\mathcal{B}$ satisfies  $($HFC$)_{k}$ with size $(M,N)$ for some $M,N\geq 2k-3$, and

\item[(iii)] $\mathbb{H}_{k}$ is weakly $(M-2k+5)$-primitive and $\mathbb{V}_{k}$ is weakly $(N-2k+5)$-primitive,
\end{enumerate}
then $\Sigma(\mathcal{B})$ has strong specification.
\end{theorem}

Theorems \ref{theorem:1.1} and \ref{theorem:1.2} are useful in verifying mixing properties. They can be used to check the results concerning strong specification and topological mixing in the literature, and can also apply to other problems. In many physical problems, edge coloring is very common. Results concerning vertex coloring can easily be extended to edge coloring and omitted here.

The rest of this paper is organized as follows.  Section 2 introduces ordering matrices of local patterns, transition matrices and connecting operators. Section 3 introduces locally corner-extendable conditions and local crisscross-extendability to study rectangle-extendability and topological mixing. Section 4 introduces invariant diagonal cycles and primitive commutative cycles of connecting operators to establish sufficient conditions for the primitivity of $\mathbb{H}_{n}$ or $\mathbb{V}_{n}$. Section 5 introduces the $k$ hole-filling condition for strong specification. The Appendix lists various mixing properties.

\numberwithin{equation}{section}

\section{Preliminaries}

\label{sec:2}

\hspace{0.5cm} This section reviews the essential aspects of the ordering matrices of local patterns and their associated transition matrices \cite{1}. It then introduces connecting operators \cite{2}.

Since the theory that was developed in this paper heavily depends on transition matrices and connecting operators, for convenience, this section presents the most important properties of transition matrices and connecting operators.

 As presented elsewhere \cite{1}, when $p\geq 2$ is fixed, the ordering matrices $\mathbf{X}_{n}$ and $\mathbf{Y}_{n}$ are introduced to arrange systematically all local patterns in $\Sigma_{2\times n}(p)$ and $\Sigma_{n\times 2}(p)$, respectively. This arrangement gives an easy recursive formulae for the transition matrices and the connecting operators, and then it gives efficient computer programming in verifying the sufficient conditions of topological mixing and strong specification. For the convenience of the readers, here we collect necessary materials from Ban and Lin \cite{1} and Ban \emph{et al.} \cite{2}.

An $n$-sequence $\overline{U}_{n}=(u_{1},u_{2},\cdots,u_{n})$ with $u_{k}\in\mathcal{S}_{p}$, $1\leq k\leq n$, is assigned a number by using the $n$-th order counting function $\psi\equiv\psi_{n}$:

\begin{equation}\label{eqn:2.1}
\psi(\overline{U}_{n})=\psi(u_{1},u_{2},\cdots,u_{n})=1+\underset{k=1}{\overset{n}{\sum}}u_{k}p^{(n-k)}.
\end{equation}
The explicit counting formula (\ref{eqn:2.1}) enables the recursive formulae that relate to $\mathbf{X}_{n}$ and $\mathbf{Y}_{n}$ to be identified.

The horizontal and vertical ordering matrices $\mathbf{X}_{2}=[x_{i_{1},j_{1}}]_{p^{2}\times p^{2}}$ and $\mathbf{Y}_{2}=[y_{i_{2},j_{2}}]_{p^{2}\times p^{2}}$ are defined by

\begin{equation}\label{eqn:2.2}
\begin{array}{ccc}
\begin{array}{c}
\psfrag{a}{$x_{i_{1},j_{1}}=$}
\psfrag{b}{{\footnotesize $u_{0,1}$}}
\psfrag{c}{{\footnotesize $u_{1,1}$}}
\psfrag{d}{{\footnotesize $u_{0,0}$}}
\psfrag{e}{{\footnotesize $u_{1,0}$}}
\psfrag{f}{}
\includegraphics[scale=1.2]{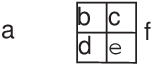}
\end{array}
&
\text{and} &
\begin{array}{c}
\psfrag{a}{$y_{i_{2},j_{2}}=$}
\psfrag{b}{{\footnotesize $u'_{0,1}$}}
\psfrag{c}{{\footnotesize $u'_{1,1}$}}
\psfrag{d}{{\footnotesize $u'_{0,0}$}}
\psfrag{e}{{\footnotesize $u'_{1,0}$}}
\psfrag{f}{,}
\includegraphics[scale=1.2]{x_ij.eps}
\end{array}
\end{array}
\end{equation}
where $u_{s,t},u'_{s,t}\in\mathcal{S}_{p}$, $0 \leq s,t\leq1$, with

\begin{equation*}
\begin{array}{ccc}
\left\{
\begin{array}{l}
i_{1}=\psi(u_{0,0},u_{0,1}) \\
j_{1}=\psi(u_{1,0},u_{1,1})
\end{array}
\right.
& \text{and} &
\left\{
\begin{array}{l}
i_{2}=\psi(u'_{0,0},u'_{1,0}) \\
j_{2}=\psi(u'_{0,1},u'_{1,1}).
\end{array}
\right.
\end{array}
\end{equation*}
For instance, when $p=2$,
\begin{equation}\label{eqn:2.2-1}
\begin{array}{ccc}
\begin{array}{c}
\psfrag{c}{$\mathbf{X}_{2}=$}
\psfrag{d}{}
\includegraphics[scale=0.475]{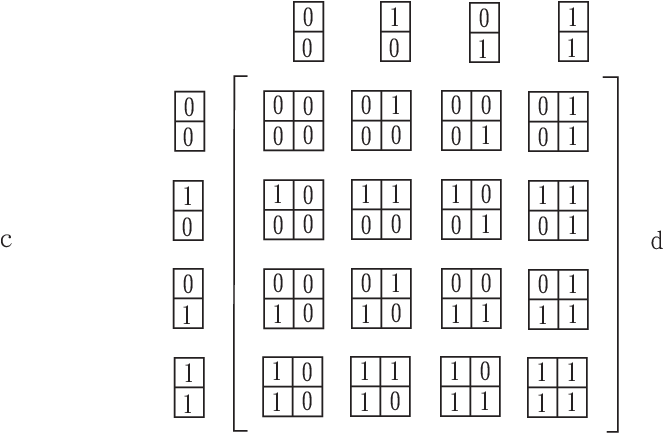}
\end{array}
& \text{and} &
\begin{array}{c}
\psfrag{c}{$\mathbf{Y}_{2}=$}
\psfrag{d}{.}
\includegraphics[scale=0.475]{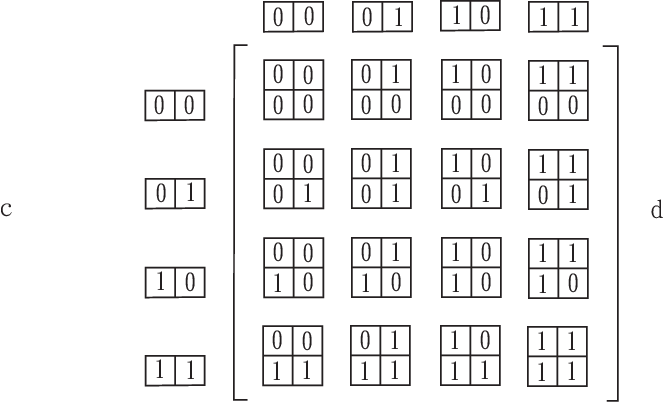}
\end{array}
\end{array}
\end{equation}

The higher-order ordering matrices $\mathbf{X}_{n}=[x_{n;i,j}]_{p^{n}\times p^{n}}$ of $\Sigma_{2\times n}(p)$, $n\geq3$, are defined recursively as

\begin{equation}\label{eqn:2.5}
\mathbf{X}_{n}=\left[X_{n;\alpha}\right]_{p\times p}=
\left[
\begin{array}{cccc}
X_{n;1} & X_{n;2} & \cdots & X_{n;p}\\
X_{n;p+1} & X_{n;p+2} & \cdots & X_{n;2p}\\
\vdots & \vdots & \ddots & \vdots \\
X_{n;p(p-1)+1} & X_{n;p(p-1)+2}  & \cdots & X_{n;p^{2}}
\end{array}
\right],
\end{equation}
where

\begin{equation}\label{eqn:2.6}
X_{n;\alpha}=
\left[
y_{\alpha, j}X_{n-1;j}
\right]_{p\times p}
\end{equation}
is a $p^{n-1}\times p^{n-1}$ matrix. Notably, the element $x_{n;i,j}$ is the $2\times n$ local pattern $U_{2\times n}=(u_{s,t})_{0\leq s\leq 1,0\leq t\leq n-1 }$ with

\begin{equation}\label{eqn:2.7}
\begin{array}{rcl}
i=\psi(u_{0,0},u_{0,1},\cdots,u_{0,n-1}) &
\text{and} &
j=\psi(u_{1,0},u_{1,1},\cdots,u_{1,n-1}).
\end{array}
\end{equation}
Similarly, the higher-order ordering matrix $\mathbf{Y}_{n}$ can be defined recursively, as above.

Given a basic set $\mathcal{B}\subset\Sigma_{2\times 2}(p)$, the horizontal and vertical transition matrices $\mathbb{H}_{2}=\mathbb{H}_{2}(B)=[h_{i,j}]_{p^{2}\times p^{2}}$ and $\mathbb{V}_{2}=\mathbb{V}_{2}(B)=[v_{i,j}]_{p^{2}\times p^{2}}$ are given by

\begin{equation}\label{eqn:2.8}
\begin{array}{ccc}
\left\{
\begin{array}{rl}
h_{i,j}=1 & \text{if }x_{i,j}\in\mathcal{B} ,\\
h_{i,j}=0 & \text{if }x_{i,j}\notin\mathcal{B},
\end{array}
\right.
&
\text{and} &
\left\{
\begin{array}{rl}
v_{i,j}=1 & \text{if }y_{i,j}\in\mathcal{B} ,\\
v_{i,j}=0 & \text{if }y_{i,j}\notin\mathcal{B}.
\end{array}
\right.
\end{array}
\end{equation}

Before the formula that relates $\mathbb{H}_{n}$ to $\mathbb{H}_{n+1}$ is presented, two kinds of products of matrices must be defined. For any two matrices $A=[a_{i,j}]$ and $B=[b_{k,l}]$, the Kronecker product (tensor product) of $A\otimes B$ is defined by

\begin{equation*}
A\otimes B=\left[a_{i,j}B\right].
\end{equation*}
Next, for any two $m\times m$ matrices $C=[c_{i,j}]$ and $D=[d_{i,j}]$, where $c_{i,j}$ and $d_{i,j}$ are numbers or matrices, the Hadamard product of $C\circ D$ is defined by

\begin{equation*}
C\circ D=\left[c_{i,j}\cdot d_{i,j}\right],
\end{equation*}
where the product $c_{i,j}\cdot d_{i,j}$ of $c_{i,j}$ and $d_{i,j}$ may be a product of numbers,
of numbers and matrices or of matrices, whenever such a product is well-defined.

According to (\ref{eqn:2.5}) and (\ref{eqn:2.6}), the higher-order transition matrices $\mathbb{H}_{n}$, $n\geq3$, can be defined as

\begin{equation}\label{eqn:2.9}
\mathbb{H}_{n}=
\left[
H_{n;\alpha}
\right]_{p\times p},
\end{equation}
where

\begin{equation}\label{eqn:2.10}
H_{n;\alpha}=
\left[
v_{\alpha,j}H_{n-1;j}
\right]_{p\times p}
\end{equation}
is a $p^{n-1}\times p^{n-1}$ zero-one matrix. Indeed, from the relation between $\mathbf{X}_{2}$ and $\mathbf{Y}_{2}$ given by (\ref{eqn:2.2}),

\begin{equation}\label{eqn:2.10-1}
H_{n;\alpha}=
\left(H_{2;\alpha}\right)_{p\times p} \circ \left[H_{n-1;j}\right]_{p\times p}.
\end{equation}

Furthermore, for any $n\geq 2$ and $q\geq 1$, $\mathbb{H}_{n+q}$ are decomposed by applying (\ref{eqn:2.9}) $q+1$ times, as follows.
For any $q\geq 1$ and $0\leq r\leq q-1$, define

\begin{equation*}
H_{n+q;\beta_{1};\beta_{2};\cdots ;\beta_{r+1}}
=
\left[
H_{n+q;\beta_{1};\beta_{2};\cdots ;\beta_{r};\alpha}
\right]_{p\times p}.
\end{equation*}
Therefore, for any $q\geq 0$, $\mathbb{H}_{n+q}$ can be represented as a $p^{q+1}\times p^{q+1}$ matrix

\begin{equation}\label{eqn:2.11}
\mathbb{H}_{n+q}\equiv\left[H_{n+q;i,j}\right]_{p^{q+1}\times p^{q+1}}=\left[H_{n+q;\beta_{1};\beta_{2};\cdots ;\beta_{q+1}}\right]_{p^{q+1}\times p^{q+1}}.
\end{equation}
In particular, when $p=2$ and $q=0$,

\begin{equation*}
\mathbb{H}_{n}=
\left[
\begin{array}{cc}
H_{n;1,1} & H_{n;1,2} \\
H_{n;2,1} & H_{n;2,2}
\end{array}
\right]_{2 \times 2}
=
\left[
\begin{array}{cc}
H_{n;1} & H_{n;2} \\
H_{n;3} & H_{n;4}
\end{array}
\right]_{2 \times 2}
;
\end{equation*}
when $p=2$ and $q=1$,
\begin{equation*}
\mathbb{H}_{n}=
\left[
\begin{array}{cccc}
H_{n;1,1} & H_{n;1,2}& H_{n;1,3}& H_{n;1,4}  \\
H_{n;2,1} & H_{n;2,2}& H_{n;2,3}& H_{n;2,4}  \\
H_{n;3,1} & H_{n;3,2}& H_{n;3,3}& H_{n;3,4}  \\
H_{n;4,1} & H_{n;4,2}& H_{n;4,3}& H_{n;4,4}
\end{array}
\right]_{2^{2} \times 2^{2}}
=
\left[
\begin{array}{cccc}
H_{n;1;1} & H_{n;1;2} & H_{n;2;1} & H_{n;2;2} \\
H_{n;1;3} & H_{n;2;4} & H_{n;2;2} & H_{n;2;4} \\
H_{n;3;1} & H_{n;3;2} & H_{n;4;1} & H_{n;4;2} \\
H_{n;3;3} & H_{n;3;4} & H_{n;4;3} & H_{n;4;4}
\end{array}
\right]_{2^{2} \times 2^{2}}
.
\end{equation*}

Now, high-order transition matrices $\mathbb{H}_{n+q}$ can be reduced to lower order transition matrices $\mathbb{H}_{n}$ as follows \cite{1}.
\begin{proposition}
\label{proposition:2.0}
For any $n\geq 2$ and $q\geq 1$,
\begin{equation}\label{eqn:2.19}
\mathbb{H}_{n+q}=\left( \mathbb{H}_{q+1} \right)_{p^{q+1}\times p^{q+1}}\circ\left(
E_{p^{q}\times p^{q}}\otimes \left[H_{n;i,j}\right]_{p\times p}
\right),
\end{equation}
where $E_{k\times k}$ is the $k\times k$ full matrix.
\end{proposition}

The formulae (\ref{eqn:2.10-1}) and (\ref{eqn:2.19}) are useful in studying rectangle-extendability, which implies that $\mathcal{B}$ is rectangle-extendable; see Theorem \ref{theorem:3.5}.

To obtain a recursive formula like that in Proposition \ref{proposition:2.0} for $\mathbb{H}_{n+q}^{m}$ to $\mathbb{H}_{n}^{m}$, $m\geq 2$ and $q\geq 1$, the connecting operator $\mathbb{C}_{m}$ must be introduced. The recursive formula is crucial in establishing sufficient conditions for the primitivity of $\mathbb{H}_{k}$ for all $k\geq 2$ by verifying the primitivity of a finite number of $\mathbb{H}_{k}$, $2\leq k\leq K$. See Theorems \ref{theorem:4.11-1} and \ref{theorem:5.3} for details.

Let $\mathbb{H}_{n}=[H_{n;i,j}]_{p\times p}$; for $m\geq 2$, the elementary pattern of $\mathbb{H}_{n}^{m}$ is

\begin{equation*}
H_{n;j_{1},j_{2}}H_{n;j_{2},j_{3}}\cdots H_{n;j_{m},j_{m+1}},
\end{equation*}
where $1\leq j_{s} \leq p$, $1\leq s\leq m+1$.
Let

\begin{equation}\label{eqn:2.20}
H_{m,n;\alpha}^{(k)}=H_{n;j_{1},j_{2}}H_{n;j_{2},j_{3}}\cdots H_{n;j_{m},j_{m+1}},
\end{equation}
where
\begin{equation*}
\begin{array}{rcl}
\alpha=\psi(j_{1}-1,j_{m+1}-1) & \text{and} & k=\psi(j_{2}-1,j_{3}-1,\cdots,j_{m}-1).
\end{array}
\end{equation*}

Therefore, for $m\geq 2$,
\begin{equation}\label{eqn:2.22}
\mathbb{H}_{n}^{m}=
\left[
H_{m,n;\alpha}
\right]_{p \times p},
\end{equation}
where

\begin{equation*}
H_{m,n;\alpha}=\underset{k=1}{\overset{p^{m-1}}{\sum}}H_{m,n;\alpha}^{(k)}.
\end{equation*}

Now, the connecting operator $\mathbb{C}_{m}=[C_{m;i,j}]$ that was introduced by \cite{2} is recalled. First, the connecting ordering matrix $\mathbf{C}_{m}=[\mathbf{C}_{m;i,j}]$ , a different arrangement for $\Sigma_{(m+1)\times 2}(p)$ from $\mathbf{Y}_{m+1}$, is introduced.
$\mathbf{C}_{m}=[\mathbf{C}_{m;i,j}]_{p^{2}\times p^{2}}$ , where $\mathbf{C}_{m;i,j}$ is a $p^{m-1}\times p^{m-1}$ matrix of local patterns, is defined as follows.

With fixed $1\leq i,j\leq p^{2}$, for $1\leq s,t \leq p^{m-1}$,

\begin{equation}\label{eqn:2.23}
\begin{array}{c}
\psfrag{a}{$\left(\mathbf{C}_{m;i,j}\right)_{s,t}=$}
\psfrag{b}{{\footnotesize $u_{0,0}$}}
\psfrag{c}{{\footnotesize $u_{1,0}$}}
\psfrag{d}{{\footnotesize$\cdots$}}
\psfrag{e}{{\footnotesize $u_{m,0}$}}
\psfrag{g}{{\footnotesize $u_{0,1}$}}
\psfrag{h}{{\footnotesize $u_{1,1}$}}
\psfrag{k}{{\footnotesize $u_{m,1}$}}
\psfrag{m}{}
\includegraphics[scale=1.3]{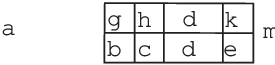}
\end{array}
\end{equation}
with $i=\psi(u_{0,0},u_{0,1})$, $j=\psi(u_{m,0},u_{m,1})$, $s=\psi(u_{1,0},u_{2,0},\cdots,u_{m-1,0})$ and $t=\psi(u_{1,1},u_{2,1},\cdots,u_{m-1,1})$.

Now, $\mathbf{C}_{m+1;i,j}$ can be obtained in terms of $\mathbf{C}_{m;k,l}$ as follows \cite{2}.

\begin{proposition}
\label{proposition:2.1}
Let $\mathbf{X}_{2}=[x_{i,j}]_{p^{2}\times p^{2}}$. For any $m\geq 2$ and $1\leq i,j\leq p^{2}$,
\begin{equation*}
\mathbf{C}_{m+1;i,j}=
\left[x_{i,\alpha}\mathbf{C}_{m;\alpha,j}
\right]_{p\times p}.
\end{equation*}
\end{proposition}

The matrix product of $\mathbf{C}_{m;i,j}$ and $\mathbf{C}_{m;j,k}$ cannot connect local patterns in the vertical direction. However, $\mathbf{S}_{m;\alpha,\beta}$ does so.
Changing the index of $\mathbf{C}_{m}=[\mathbf{C}_{m;i,j}]_{p^{2}\times p^{2}}$ enables the ordering matrix $\mathbf{S}_{m}=[\mathbf{S}_{m;\alpha,\beta}]_{p^{2}\times p^{2}}$ to be defined as

\begin{equation}\label{eqn:2.24}
\mathbf{S}_{m;\alpha,\beta}=\mathbf{C}_{m;\psi(\alpha_{1},\beta_{1}),\psi(\alpha_{2},\beta_{2})},
\end{equation}
where $\alpha_{k},\beta_{k}\in\mathcal{S}_{p}$, $1\leq k\leq 2$, satisfying $\alpha=\psi(\alpha_{1},\alpha_{2})$ and $\beta=\psi(\beta_{1},\beta_{2})$.
In particular, for $p=2$,
\begin{equation*}
\mathbf{C}_{m}=
\left[
\begin{array}{cccc}
\mathbf{C}_{m;1,1} & \mathbf{C}_{m;1,2} & \mathbf{C}_{m;1,3} & \mathbf{C}_{m;1,4} \\
\mathbf{C}_{m;2,1} & \mathbf{C}_{m;2,2} & \mathbf{C}_{m;3,3} & \mathbf{C}_{m;4,4} \\
\mathbf{C}_{m;3,1} & \mathbf{C}_{m;3,2} & \mathbf{C}_{m;3,3} & \mathbf{C}_{m;3,4} \\
\mathbf{C}_{m;4,1} & \mathbf{C}_{m;4,2} & \mathbf{C}_{m;4,3} & \mathbf{C}_{m;4,4}
\end{array}
\right]
=
\left[
\begin{array}{cccc}
\mathbf{S}_{m;1,1} & \mathbf{S}_{m;1,2} & \mathbf{S}_{m;2,1} & \mathbf{S}_{m;2,2} \\
\mathbf{S}_{m;1,3} & \mathbf{S}_{m;1,4} & \mathbf{S}_{m;2,3} & \mathbf{S}_{m;2,4} \\
\mathbf{S}_{m;3,1} & \mathbf{S}_{m;3,2} & \mathbf{S}_{m;4,1} & \mathbf{S}_{m;4,2} \\
\mathbf{S}_{m;3,3} & \mathbf{S}_{m;3,4} & \mathbf{S}_{m;4,3} & \mathbf{S}_{m;4,4}
\end{array}
\right].
\end{equation*}
Indeed, for $1\leq s,t \leq p^{m-1}$,

\begin{equation}\label{eqn:2.25}
\begin{array}{c}
\psfrag{a}{$(\mathbf{S}_{m;\alpha,\beta})_{s,t}=$}
\psfrag{b}{{\footnotesize $u_{0,0}$}}
\psfrag{c}{{\footnotesize $u_{1,0}$}}
\psfrag{d}{{\footnotesize$\cdots$}}
\psfrag{e}{{\footnotesize $u_{m,0}$}}
\psfrag{g}{{\footnotesize $u_{0,1}$}}
\psfrag{h}{{\footnotesize $u_{1,1}$}}
\psfrag{k}{{\footnotesize $u_{m,1}$}}
\psfrag{m}{}
\includegraphics[scale=1.4]{C_m_ij_st.eps}
\end{array}
\end{equation}
with $\alpha=\psi(u_{0,0},u_{m,0})$, $\beta=\psi(u_{0,1},u_{m,1})$, $s=\psi(u_{1,0},u_{2,0},\cdots,u_{m-1,0})$ and $t=\psi(u_{1,1},u_{2,1},\cdots,u_{m-1,1})$. From (\ref{eqn:2.25}), the matrix product of $\mathbf{S}_{m;\alpha,\beta}$ and
$\mathbf{S}_{m;\beta,\gamma}$ represents the vertical connection of the patterns on $\mathbb{Z}_{(m+1)\times 2}$.

Now, given $\mathcal{B}\subset\Sigma_{2\times 2}(p)$, for $m\geq 2$, the connecting operator $\mathbb{C}_{m}=[C_{m;i,j}]_{p^{2}\times p^{2}}$ of  $\mathbf{C}_{m}=[\mathbf{C}_{m;i,j}]_{p^{2}\times p^{2}}$ is defined as follows. For $1\leq s,t \leq p^{m-1}$,

\begin{equation*}
\left\{
\begin{array}{rl}
(C_{m;i,j})_{s,t}=1 &   \text{if }(\mathbf{C}_{m;i,j})_{s,t} \text{ is }\mathcal{B}\text{-admissible}, \\
(C_{m;i,j})_{s,t}=0 &  \text{otherwise.}
\end{array}
\right.
\end{equation*}

In the following, $\mathbb{C}_{2}$ can be obtained explicitly. Let $\mathbb{H}_{2}=[h_{i,j}]_{p^{2}\times p^{2}}$.
Then, for $1\leq i,j\leq p^{2}$,

{\footnotesize
\begin{equation}\label{eqn:2.26}
\begin{array}{rl}
 & C_{2;i,j} \\
 & \\
 = &
\left[
\begin{array}{cccc}
h_{i,1} & h_{i,2} & \cdots & h_{i,p} \\
h_{i,p+1} & h_{i,p+2} &\cdots & h_{i,2p} \\
\vdots& \vdots & \vdots & \vdots \\
h_{i,(p-1)p+1} & h_{i,(p-1)p+2} & \cdots & h_{i,p^{2}}
\end{array}
\right] \circ
\left[
\begin{array}{cccc}
h_{1,j} & h_{2,j} & \cdots & h_{p,j} \\
h_{p+1,j} & h_{p+2,j} &\cdots & h_{2p,j} \\
\vdots& \vdots & \vdots & \vdots \\
h_{(p-1)p+1,j} & h_{(p-1)p+2,j} & \cdots & h_{p^{2},j}
\end{array}
\right]
\end{array}
\end{equation}
}
is a $p\times p$ zero-one matrix. By Proposition \ref{proposition:2.1}, the connecting operator $\mathbb{C}_{m+1}$ can also be obtained from $\mathbb{C}_{m}$. For $m\geq 2$, $\mathbb{C}_{m+1}=[C_{m+1;i,j}]_{p^{2}\times p^{2}}$ satisfies

\begin{equation}\label{eqn:2.27}
C_{m+1;i,j}=
\left[
h_{i,\alpha}C_{m;\alpha,j}
\right]_{p\times p}.
\end{equation}

From (\ref{eqn:2.24}), $\mathbb{S}_{m}=[S_{m;\alpha,\beta}]_{p^{2}\times p^{2}}$ is defined by

\begin{equation}\label{eqn:2.27-1}
S_{m;\alpha,\beta}=C_{m;\psi(\alpha_{1},\beta_{1}),\psi(\alpha_{2},\beta_{2})},
\end{equation}
where $0\leq\alpha_{1},\alpha_{2},\beta_{1},\beta_{2}\leq p-1$ such that $\alpha=\psi(\alpha_{1},\alpha_{2})$ and $\beta=\psi(\beta_{1},\beta_{2})$.

For example, consider the Golden Mean shift,
\begin{equation*}
\mathbb{H}_{2}=\mathbb{V}_{2}=\left[
\begin{array}{cccc}
1 & 1 & 1 & 0 \\
1 & 0 & 1 & 0 \\
1 & 1 & 0 & 0 \\
0 & 0 & 0 & 0
\end{array}
\right].
\end{equation*}
By (\ref{eqn:2.28}), it can be verified that

\begin{equation*}
\begin{array}{rl}
\mathbb{C}_{2}= &
\left[
\begin{array}{cccc}
C_{2;1,1} & C_{2;1,2} & C_{2;1,3} & C_{2;1,4} \\
C_{2;2,1} & C_{2;2,2} & C_{2;2,3} & C_{2;2,4} \\
C_{2;3,1} & C_{2;3,2} & C_{2;3,3} & C_{2;3,4} \\
C_{2;4,1} & C_{2;4,2} & C_{2;4,3} & C_{2;4,4}
\end{array}
\right]
\\
& \\
= &
{\scriptsize
\left[
\begin{array}{cccc}
\left[\begin{array}{cc} 1& 1 \\ 1& 0 \end{array}\right] &\left[\begin{array}{cc} 1& 0 \\ 1& 0 \end{array}\right] & \left[\begin{array}{cc} 1& 1 \\ 0& 0 \end{array}\right] & \left[\begin{array}{cc} 0& 0 \\0& 0 \end{array}\right] \\
& & & \\
\left[\begin{array}{cc} 1& 0 \\ 1& 0 \end{array}\right] &\left[\begin{array}{cc} 1& 0 \\ 1& 0 \end{array}\right] & \left[\begin{array}{cc} 1& 0 \\ 0& 0 \end{array}\right] & \left[\begin{array}{cc} 0& 0 \\ 0& 0 \end{array}\right] \\
& & & \\
\left[\begin{array}{cc} 1& 1 \\ 0& 0 \end{array}\right] &\left[\begin{array}{cc} 1& 0 \\ 0& 0 \end{array}\right] & \left[\begin{array}{cc} 1& 1 \\ 0& 0 \end{array}\right] & \left[\begin{array}{cc} 0& 0 \\ 0& 0 \end{array}\right] \\
& & & \\
\left[\begin{array}{cc} 0& 0 \\ 0& 0 \end{array}\right] &\left[\begin{array}{cc} 0& 0 \\ 0& 0 \end{array}\right] & \left[\begin{array}{cc} 0& 0 \\ 0& 0 \end{array}\right] & \left[\begin{array}{cc} 0& 0 \\ 0& 0 \end{array}\right]
\end{array}
\right]. }
\end{array}
\end{equation*}
Moreover, by (\ref{eqn:2.27-1}),
\begin{equation*}
\mathbb{S}_{2}= \left[
\begin{array}{cccc}
S_{2;1,1} & S_{2;1,2} & S_{2;1,3} & S_{2;1,4} \\
S_{2;2,1} & S_{2;2,2} & S_{2;2,3} & S_{2;2,4} \\
S_{2;3,1} & S_{2;3,2} & S_{2;3,3} & S_{2;3,4} \\
S_{2;4,1} & S_{2;4,2} & S_{2;4,3} & S_{2;4,4}
\end{array}
\right]
=\left[
\begin{array}{cccc}
C_{2;1,1} & C_{2;1,2} & C_{2;2,1} & C_{2;2,2} \\
C_{2;1,3} & C_{2;1,4} & C_{2;2,3} & C_{2;2,4} \\
C_{2;3,1} & C_{2;3,2} & C_{2;4,1} & C_{2;4,2} \\
C_{2;3,3} & C_{2;3,4} & C_{2;4,3} & C_{2;4,4}
\end{array}
\right].
\end{equation*}

Now, the relation between $\mathbb{H}_{n+1}^{m}$ and $\mathbb{H}_{n}^{m}$ is elucidated as follows. Since the sizes of $H_{m,n+1;\alpha}^{(k)}$ and $H_{m,n;\beta}^{(l)}$ are different, the elementary pattern $H_{m,n+1;\alpha}^{(k)}$ can be reduced further as follows.

Let

\begin{equation}\label{eqn:2.28}
H_{m,n+1;\alpha}^{(k)}=
\left[H_{m,n+1;\alpha;\beta}^{(k)}
\right]_{p\times p}.
\end{equation}

In the following theorem, $H^{(k)}_{m,n+1;\alpha;\beta}$ is obtained as the product of $S_{m;\alpha,\beta}$ and $H^{(l)}_{m,n;\beta}$ \cite{2}, so $S_{m;\alpha,\beta}$ reduces $\mathbb{H}_{n+1}^{m}$ to $\mathbb{H}_{n}^{m}$.
\begin{proposition}
\label{proposition:2.2}
For any $m,n\geq 2$,

\begin{equation}\label{eqn:2.31}
H^{(k)}_{m,n+1;\alpha;\beta}= \underset{l=1}{\overset{p^{m-1}}{\sum}}\left( S_{m;\alpha,\beta}\right)_{k,l}H^{(l)}_{m,n;\beta}.
\end{equation}
Furthermore, for $n=1$, let

\begin{equation}\label{eqn:2.32}
H_{m,2;\alpha}^{(k)}=
\left[ H_{m,2;\alpha;\beta}^{(k)}
\right]_{p \times p},
\end{equation}
then
\begin{equation}\label{eqn:2.33}
H_{m,2;\alpha;\beta}^{(k)}=\underset{l=1}{\overset{p^{m-1}}{\sum}}\left(S_{m;\alpha,\beta}\right)_{k,l}.
\end{equation}
\end{proposition}

Furthermore, for $q\geq 2$, $q$-many $S_{m;\alpha,\beta}$ can reduce $\mathbb{H}_{n+q}^{m}$ to $\mathbb{H}_{n}^{m}$ as follows.
For any positive integer $q\geq 2$, the elementary patterns of $\mathbb{H}_{n+q}^{m}$ can be decomposed by applying (\ref{eqn:2.28}) $q$ times. Indeed,
for $q\geq 2$ and $1\leq r\leq q-1$, define

\begin{equation*}
H_{m,n+q;\beta_{1};\beta_{2};\cdots ;\beta_{r+1}}^{(k)}
 =
\left[
 H_{m,n+q;\beta_{1};\beta_{2};\cdots ;\beta_{r+1};\alpha}^{(k)}
\right]_{p\times p}.
\end{equation*}
Therefore, for any $q\geq 1$, $\mathbb{H}_{n+q}^{m}$ can be represented as a $p^{q+1}\times p^{q+1}$ matrix

\begin{equation}\label{eqn:2.34}
\mathbb{H}_{n+q}^{m}\equiv \left[H_{m,n+q;i,j}\right]_{p^{q+1}\times p^{q+1}}=\left[H_{m,n+q;\beta_{1};\beta_{2};\cdots ;\beta_{q+1}}\right]_{p^{q+1}\times p^{q+1}}
\end{equation}
where

\begin{equation*}
H_{m,n+q;\beta_{1};\beta_{2};\cdots ;\beta_{q+1}}=\underset{k=1}{\overset{p^{m-1}}{\sum}}H_{m,n+q;\beta_{1};\beta_{2};\cdots ;\beta_{q+1}}^{(k)}
\end{equation*}
is a $p^{n-1}\times p^{n-1}$ matrix.
%

 As in Proposition \ref{proposition:2.2}, the elementary patterns of $\mathbb{H}_{n+q}^{m}$ can be expressed as the product of $q$-many $S_{m;\alpha,\beta}$ and the elementary patterns of $\mathbb{H}_{n}^{m}$ \cite{2}.

\begin{proposition}
\label{proposition:2.3}
For any $m,n\geq 2$ and $q\geq  1$,

\begin{equation}\label{eqn:2.35}
H_{m,n+q;\beta_{1};\beta_{2};\cdots ;\beta_{q+1}}^{(k)}=\underset{l=1}{\overset{p^{m-1}}{\sum}}
(S_{m;\beta_{1},\beta_{2}}S_{m;\beta_{2},\beta_{3}}\cdots S_{m;\beta_{q},\beta_{q+1}})_{k,l}
H_{m,n;\beta_{q+1}}^{(l)}.
\end{equation}
where $ 1 \leq \beta_{i}\leq p^{2}$, $1\leq i\leq q+1$. Moreover,

\begin{equation}\label{eqn:2.36}
H_{m,n+q;\beta_{1};\beta_{2};\cdots ;\beta_{q+1}}=\underset{k,l=1}{\overset{p^{m-1}}{\sum}}
(S_{m;\beta_{1},\beta_{2}}S_{m;\beta_{2},\beta_{3}}\cdots S_{m;\beta_{q},\beta_{q+1}})_{k,l}
H_{m,n;\beta_{q+1}}^{(l)}.
\end{equation}
\end{proposition}

Similarly, for $\mathbb{V}_{2}$, the connecting operators are denoted by $\mathbb{U}_{m}=[U_{m;i,j}]$ (corresponding to $\mathbb{C}_{m}=[C_{m;i,j}]$ for $\mathbb{H}_{2}$) and $\mathbb{W}_{m}=[W_{m;\alpha,\beta}]$ (corresponding to $\mathbb{S}_{m}=[S_{m;\alpha,\beta}]$ for $\mathbb{H}_{2}$). The arguments that hold for $\mathbb{H}_{n}$ also hold for $\mathbb{V}_{n}$.

In the study of both topological mixing and strong specification, the transition matrices $\mathbb{H}_{n}$ and $\mathbb{V}_{n}$ and the connecting operator $\mathbb{S}_{m}$ or $\mathbb{C}_{m}$ are extensively used. Indeed, invariant diagonal cycles, primitive commutative cycles and (HFC$)_{k}$ can be expressed in terms of transition matrices and connecting operators as the finitely sufficient conditions; see Definitions 4.1 and 4.6 and Theorem 5.2.
All cases with certain extendability conditions except strong specification can be verified using transition matrices and connecting operators. Table 2.1 lists the related theorems.

\begin{equation*}
\psfrag{a}{Mixing properties}
\psfrag{b}{Results }
\psfrag{c}{Strong specification}
\psfrag{d}{Uniform filling property }
\psfrag{e}{Corner gluing}
\psfrag{f}{Block gluing}
\psfrag{g}{Topological mixing}
\psfrag{h}{{\small Expressions in $\mathbb{H}$ and $\mathbb{S}$}}
\psfrag{j}{Sufficient conditions in}
\psfrag{z}{{\small finitely expressions}}
\psfrag{k}{Yes}
\psfrag{l}{Yes}
\psfrag{m}{Yes}
\psfrag{n}{Theorem \ref{theorem:1.1}}
\psfrag{o}{Theorem \ref{theorem:1.2}}
\psfrag{p}{Theorem \ref{theorem:5.7}}
\includegraphics[scale=1.07]{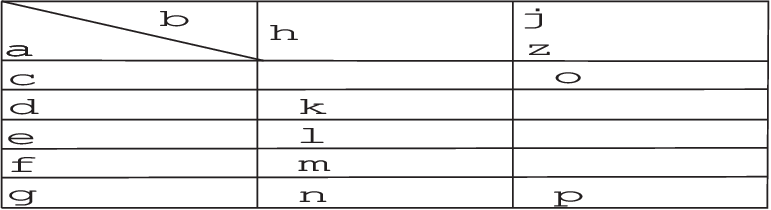}
\end{equation*}

\begin{equation*}
\text{Table 2.1.}
\end{equation*}

\numberwithin{equation}{section}

\section{Extendabilities and topological mixing}

\label{sec:3}
\hspace{0.5cm} This section investigates the extendabilities and the relationship with the topological mixing of $\Sigma(\mathcal{B})$.

First, the main idea of finding sufficient conditions for topological mixing is stated. Given two patterns $U_{1}$ and $U_{2}$ defined on $R_{1}$ and $R_{2}+\mathbf{v}$, respectively, in general, $R_{1}$ and $R_{2}+\mathbf{v}$ are not located in a horizontal line or a vertical line. Typically, the gluing process comprises three steps; an example is presented in Fig. 3.1. For clarity, in Fig 3.1, the patterns, $U$, are presented and the underlying lattices, $R$, are omitted.

\begin{enumerate}
\item[Step (1):] Extend $U_{2}$ to $\widetilde{U}_{2}$ such that $U_{1}$ can connect horizontally to $\widetilde{U}_{2}$. The combined pattern is an $L$-shaped pattern $U_{1}\bigcup\widetilde{U}_{1}\bigcup U_{2}\bigcup\widetilde{U}_{2}$.

\item[Step (2):] Extend the $L$-shaped pattern to a rectangular block, $U_{1}\bigcup\widetilde{U}_{1}\bigcup U_{2}\bigcup\widetilde{U}_{2}\bigcup U_{3}$.

\item[Step (3):] Extend the rectangular block to a global pattern on $\mathbb{Z}^{2}$.
\end{enumerate}

\begin{equation*}
\psfrag{a}{{\footnotesize $U_{1}$}}
\psfrag{b}{{\footnotesize$U_{2}$}}
\psfrag{c}{{\footnotesize$\widetilde{U}_{2}$}}
\psfrag{d}{{\footnotesize $(1)$}}
\psfrag{e}{{\footnotesize$\widetilde{U}_{1}$}}
\psfrag{f}{\hspace{-0.5cm}{\footnotesize$U_{3}$}\hspace{0.2cm}$(2)$}
\psfrag{g}{ $(3)$}
\includegraphics[scale=0.9]{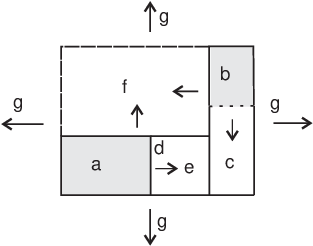}
\end{equation*}
\begin{equation*}
\text{Figure 3.1.}
\end{equation*}

To ensure that all processes can be executed, the following sufficient conditions are proposed to be applied in each step:

In Step (1), since $U_{2}$ is part of a global pattern, $U_{2}$ can be extended to $\widetilde{U}_{2}$.
The introduction of the primitivity of horizontal transition matrices $\mathbb{H}_{n}$ and vertical transition matrices $\mathbb{V}_{n}$, for each $n\geq 2$, ensures $U_{1}$ can connect to $\widetilde{U}_{2}$.

In Step (2), corner-extendability is introduced to enable every admissible $L$-shaped pattern to be extended to be a rectangular block as $U_{1}\bigcup\widetilde{U}_{1}\bigcup U_{2}\bigcup\widetilde{U}_{2}\bigcup U_{3}$.

In Step (3), rectangle-extendability is introduced to extend every rectangular block to form a global pattern on $\mathbb{Z}^{2}$; see Theorem \ref{theorem:3.5}.

Notably, extending $U_{2}$ to $\widetilde{U}_{2}$ and $U_{1}$ to $\widetilde{U}_{1}$ simultaneously demands a stronger sufficient condition to connect $\widetilde{U}_{1}$ and $\widetilde{U}_{2}$, meaning that there exists a constant $M\geq1$ such that $\mathbb{H}_{n}$ and $\mathbb{V}_{n}$ are $M$-primitive for all $n\geq 2$.  Actually, this condition is like block gluing (see Definition A.1 (iv)). Therefore, the separate execution of Steps (1) and (2) weakens the sufficient conditions for topological mixing.

After the primitivity of $\mathbb{H}_{n}$ and $\mathbb{V}_{n}$ has been established, the locally corner-extendable conditions $C(1)\sim C(4)$ and locally crisscross-extendability are introduced to extend the L-shaped pattern and the rectangular pattern into a global pattern, and then to establish that $\Sigma(\mathcal{B})$ is topologically mixing.

The importance of the corners of a finite lattice has been noticed \cite{29-1, 42}. Johnson \emph{et al.} \cite{29-1} introduced the concept of corner gluing in a study of factors of higher-dimensional shifts of finite type. Similarly, to study rectangle-extendability and topological mixing, the corners of a rectangular lattice must be studied closely. Indeed, let the $L$-shaped lattices $\mathbb{L}_{1}=\mathbb{Z}_{3\times 3} \setminus \{(2,2)\}$, $\mathbb{L}_{2}=\mathbb{Z}_{3\times 3} \setminus \{(0,2)\}$,  $\mathbb{L}_{3}=\mathbb{Z}_{3\times 3} \setminus \{(0,0)\}$ and $\mathbb{L}_{4}=\mathbb{Z}_{3\times 3} \setminus \{(2,0)\}$; accordingly,
\begin{equation*}
\begin{array}{cccc}
\psfrag{a}{{\footnotesize $\mathbb{L}_{1}= $}}
\includegraphics[scale=0.7]{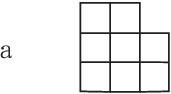},  &
\psfrag{a}{{\footnotesize$\mathbb{L}_{2}=$}}
\includegraphics[scale=0.7]{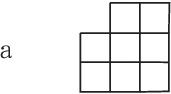},
  &
\psfrag{a}{\footnotesize{$\mathbb{L}_{3}=$}}
\includegraphics[scale=0.7]{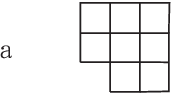},
  &
\psfrag{a}{\footnotesize{$\mathbb{L}_{4}=$}}
\includegraphics[scale=0.7]{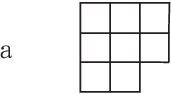}
\end{array}.
\end{equation*}
\begin{equation*}
\text{Figure 3.2.}
\end{equation*}
For $1\leq i\leq 4$, a given $\mathcal{B}\subset\Sigma_{2\times 2}(p)$ satisfies the locally corner-extendable condition $C(i)$ if
for any $U\in\Sigma_{\mathbb{L}_{i}}(\mathcal{B})$, there exists $U'\in\Sigma_{3\times 3}(\mathcal{B})$ such that $U'\mid_{\mathbb{L}_{i}}=U$.

The crisscross lattice $\mathbb{Z}_{c}$ is defined by

\begin{equation}\label{eqn:3.1}
\mathbb{Z}_{c}=\underset{0\leq |i|+|j|\leq 1}{\bigcup}\mathbb{Z}_{2\times 2}((i,j)).
\end{equation}
Indeed,
\begin{equation*}
\psfrag{a}{{\tiny $O$}}
\psfrag{f}{$\mathbb{Z}_{c}=$}
\psfrag{g}{,}
\includegraphics[scale=0.8]{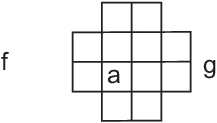}
\end{equation*}
\begin{equation*}
\text{Figure 3.3.}
\end{equation*}
where $O=(0,0)$ is the origin of $\mathbb{Z}^{2}$.
For $\mathcal{B}\subset\Sigma_{2\times 2}(p)$, let
\begin{equation*}
\Sigma_{c}(\mathcal{B})=\Sigma_{\mathbb{Z}_{c} }(\mathcal{B}).
\end{equation*}
For $\mathcal{B}\subset\Sigma_{2\times 2}(p)$, $\mathcal{B}$ satisfies local crisscross-extendability if each $B\in \mathcal{B}$, there exists $U_{c}\in\Sigma_{c}(\mathcal{B})$ with $U_{c}\mid_{\mathbb{Z}_{2\times 2}}=B$.

Clearly, primitivity of $\mathbb{H}_{n}$ and $\mathbb{V}_{n}$ for all $n\geq 2$ can be interpreted as topological mixing in horizontal and vertical directions, respectively. In general, the lattices $R_{1}$ and $R_{2}+\mathbf{v}$ are not located along horizontal or vertical lines. Accordingly, mixing in directions other than horizontal and vertical must be studied. In so doing, locally corner-extendable conditions and local crisscross-extendability are useful.

First, the rectangle-extendability of $\mathcal{B}$ is defined as follows.

\begin{definition}
\label{definition:3.1}
For $\mathcal{B}\subset\Sigma_{2\times 2}(p)$, $\mathcal{B}$ is called rectangle-extendable if for every rectangular block $U_{m\times n}\in\Sigma_{m\times n}(\mathcal{B})$, $m,n\geq 2$, there exists $W\in\Sigma(\mathcal{B})$ such that $W\mid_{\mathbb{Z}_{m\times n}}=U_{m\times n}$.
\end{definition}
In general, due to the undecidability of two-dimensional shifts of finite type, the rectangle-extendability is not finitely checkable.
Notably, if $\mathcal{B}$ is rectangle-extendable, then $\Sigma(\mathcal{B})\neq \emptyset$. The converse is not true in general.

Whether or not $\mathbb{H}_{2}(\mathcal{B})$ or $\mathbb{V}_{2}(\mathcal{B})$ contains a zero row or a zero column has a very large impact in studying mixing problems. First, the case in which a matrix contains no zero row or zero column is considered.

\begin{definition}
\label{definition:3.3}
A matrix $A=[a_{i,j}]_{n\times n}$ is non-compressible if it contains no zero row and no zero column.
For $n\geq 2$, an $\mathbb{H}_{n}$ (or $\mathbb{V}_{n}$) is non-degenerated if $H_{n;\alpha}$ (or $V_{n;\alpha}$) is non-compressible for all $1\leq\alpha\leq p^{2}$.
\end{definition}

First, consider the case of $\mathcal{B}\subset \Sigma_{2\times 2}(p)$ when $\mathbb{H}_{2}(\mathcal{B})$ and  $\mathbb{V}_{2}(\mathcal{B})$ are non-degenerated.
Clearly, if both $A$ and $B$ are non-negative and non-compressible matrices, then $AB$ is non-compressible. In the following, the recursive formula from high-order transition matrices $\mathbb{H}_{n}$ and $\mathbb{V}_{n}$ to $\mathbb{H}_{2}$ and $\mathbb{V}_{2}$ is proven to ensure that $\mathbb{H}_{n}$ (or $\mathbb{V}_{n}$) is non-degenerated for $n\geq 3$ if it holds for $\mathbb{H}_{2}$ (or $\mathbb{V}_{2}$).

\begin{theorem}
\label{theorem:3.4-0}
If $\mathbb{H}_{2}$ is non-degenerated, then $\mathbb{H}_{n}$ is non-degenerated for all $n\geq 3$.  Moreover, $H_{m,n;\alpha}^{(k)}$ are also non-compressible for $m, n\geq 2$, $1\leq \alpha\leq p^{2}$ and $1\leq k\leq p^{m-1}$.
\end{theorem}

\textit{Proof.}
Since $\mathbb{H}_{2}$ is non-degenerated, from Definition  \ref{definition:3.3}, $H_{2;\alpha}$ are non-compressible for $1\leq \alpha\leq p^{2}$. For $n\geq 3$ and $1\leq \alpha\leq p^{2}$, (\ref{eqn:2.10-1}) implies $H_{n;\alpha}$ are non-compressible. Then, $\mathbb{H}_{n}$ is non-degenerated for all $n\geq 3$.  From (\ref{eqn:2.20}), that $H_{m,n;\alpha}^{(k)}$ is non-compressible follows. \hspace{0.5cm} $\square$
\medbreak

The following theorem provides sufficient conditions on $\mathbb{H}_{n}$ and $\mathbb{V}_{n}$ for rectangle-extendability.

\begin{theorem}
\label{theorem:3.5}
Given $\mathcal{B}\subset\Sigma_{2\times 2}(p)$, if $\mathbb{H}_{n}(\mathcal{B})$ and $\mathbb{V}_{n}(\mathcal{B})$ are non-compressible for all $n\geq 2$, then
$\mathcal{B}$ is rectangle-extendable. Furthermore, if $\mathbb{H}_{2}(\mathcal{B})$ and $\mathbb{V}_{2}(\mathcal{B})$ are non-degenerated, then $\mathcal{B}$ is rectangle-extendable. In particular, $\Sigma(\mathcal{B})\neq \emptyset$.
\end{theorem}

\textit{Proof.}
Given $U\in \Sigma_{m\times n}(\mathcal{B})$, if $\mathbb{H}_{n}$ is non-compressible, then $U$ can be extended in both positive and negative horizontal directions to form an $(m+2)\times n$ $\mathcal{B}$-admissible pattern $U_{1}$. Similarly, the fact that $\mathbb{V}_{m+2}$ is non-compressible implies that $U_{1}$ can be extended to an $(m+2)\times (n+2)$ $\mathcal{B}$-admissible pattern $U_{2}$. Repeating this process extends $U$ to a global pattern in $\Sigma(\mathcal{B})$. Therefore, $\mathcal{B}$ is rectangle-extendable.
Moreover, if $\mathbb{H}_{2}(\mathcal{B})$ and $\mathbb{V}_{2}(\mathcal{B})$ are non-degenerated, then the result follows from Theorem \ref{theorem:3.4-0}. \hspace{0.5cm} $\square$
\medbreak

The non-degeneracy of $\mathbb{H}_{2}(\mathcal{B})$ and $\mathbb{V}_{2}(\mathcal{B})$ implies three of the locally corner-extendable conditions, as follows.

\begin{theorem}
\label{theorem:3.6}
Given $\mathcal{B}\subset \Sigma_{2\times 2}(p)$,
if $\mathbb{H}_{2}(\mathcal{B})$ and $\mathbb{V}_{2}(\mathcal{B})$ are non-degenerated, then
$\mathcal{B}$ satisfies $C(1)$, $C(2)$ and $C(4)$.
\end{theorem}

\textit{Proof.}
Since $\mathbb{H}_{2}(\mathcal{B})$ is non-degenerated, from (\ref{eqn:2.2}), for any $u_{0,0},u_{1,0},u_{0,1},u_{1,1}\in\mathcal{S}_{p}$, there exist $a,b\in\mathcal{S}_{p}$ such that

\begin{equation*}
\begin{array}{ccc}
\psfrag{c}{$a$}
\psfrag{b}{{\footnotesize $u_{0,1}$}}
\psfrag{d}{{\footnotesize $u_{0,0}$}}
\psfrag{e}{{\footnotesize $u_{1,0}$}}
\includegraphics[scale=1.2]{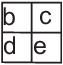} &\hspace{0.5cm}
\text{and} & \hspace{0.5cm}
\psfrag{b}{$b$}
\psfrag{c}{{\footnotesize $u_{1,1}$}}
\psfrag{d}{{\footnotesize $u_{0,0}$}}
\psfrag{e}{{\footnotesize $u_{1,0}$}}
\includegraphics[scale=1.2]{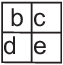}
\end{array}
\end{equation*}
\begin{equation*}
\text{Figure 3.4.}
\end{equation*}
are in $\mathcal{B}$, which implies that conditions $C(1)$ and $C(2)$ are satisfied.
Similarly, that $\mathbb{V}_{2}(\mathcal{B})$ is non-degenerated implies that $\mathcal{B}$ satisfies conditions $C(1)$ and $C(4)$.

The proof is complete. \hspace{0.5cm} $\square$
\medbreak

Now, the fact that $\Sigma(\mathcal{B})$ is topologically mixing follows from the non-degeneracy of $\mathbb{H}_{2}(\mathcal{B})$ and $\mathbb{V}_{2}(\mathcal{B})$ and the primitivity of $\mathbb{H}_{n}$ and $\mathbb{V}_{n}$, $n\geq 2$.

\begin{theorem}
\label{theorem:3.7}
Given $\mathcal{B}\subset \Sigma_{2\times 2}(p)$, if $\mathbb{H}_{2}(\mathcal{B})$ and $\mathbb{V}_{2}(\mathcal{B})$ are non-degenerated, then the following statements are equivalent.

\item[(i)] $\mathbb{H}_{n}(\mathcal{B})$ and $\mathbb{V}_{n}(\mathcal{B})$ are primitive for all $n\geq 2$.

\item[(ii)]$\Sigma(\mathcal{B})$ is topologically mixing.
\end{theorem}

\textit{Proof.}
(i)$\Rightarrow$(ii). Let $R_{1}$ and $R_{2}$ be finite sublattices of $\mathbb{Z}^{2}$. Then, there exist $N\geq 2$ and $({i_{1},j_{1}})$,$(i_{2},j_{2})\in\mathbb{Z}^{2}$ such that $R_{l}\subset\mathbb{Z}_{N\times N}((i_{l},j_{l}))$, $l=1,2$. From (i), there exists $K\geq1$ such that $\mathbb{H}_{N}^{K}(\mathcal{B})>0$ and $\mathbb{V}_{N}^{K}(\mathcal{B})>0$.

Then, consider $M=M(R_{1},R_{2})=\sqrt{2}(2N+K-2)$. Let $\mathbf{v}=(v_{1},v_{2})\in\mathbb{Z}^{2}$ with $d(R_{1},R_{2}+\mathbf{v})\geq M$ and any two allowable patterns $U_{1}\in\Pi_{R_{1}}(\Sigma(\mathcal{B}))$ and $U_{2}\in\Pi_{R_{2}+\mathbf{v}}(\Sigma(\mathcal{B}))$.
Clearly, $U_{1}$ and $U_{2}$ can be extended as $U_{1}'$ on $\mathbb{Z}_{N\times N}((i_{1},j_{1}))$ and $U_{2}'$ on $\mathbb{Z}_{N\times N}((i_{2}+v_{1},j_{2}+v_{2}))$ using the local patterns in $\mathcal{B}$, respectively.

Proving that $U_{1}'$ and  $U_{2}'$ can be connected to form the L-shaped pattern $U_{L}$ using the local patterns in $\mathcal{B}$, as follows, is not difficult.
\begin{equation*}
\begin{array}{lcr}
\psfrag{a}{{\footnotesize$N$}}
\psfrag{b}{$U_{L}$}
\includegraphics[scale=0.65]{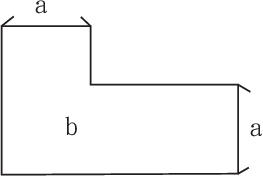} & \hspace{0.5cm}\text{or} &\hspace{0.5cm}
\psfrag{a}{{\footnotesize$N$}}
\psfrag{b}{$U_{L}$}
\includegraphics[scale=0.65]{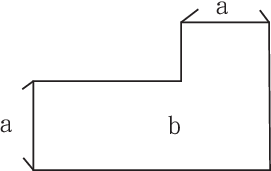}
\end{array}
\end{equation*}
\begin{equation*}
\text{Figure 3.5.}
\end{equation*}
Notably, the L-shaped lattices may degenerate into rectangular lattices.

Since $\mathbb{H}_{2}(\mathcal{B})$ and $\mathbb{V}_{2}(\mathcal{B})$ are non-degenerated, by Theorem \ref{theorem:3.6}, $\mathcal{B}$ satisfies conditions $C(1)$ and $C(2)$. Then, $U_{L}$ can be extended as $U_{r}$ on the rectangular lattice by using the local patterns in $\mathcal{B}$, which is obtained by filling the corner of the L-shaped lattices.

From Theorem \ref{theorem:3.4-0}, $\mathcal{B}$ is rectangle-extendable.
Then, $U_{r}$ can be extended as $W\in\Sigma(\mathcal{B})$ with $\Pi_{R_{1}}(W)=U_{1}$ and $\Pi_{R_{2}+\mathbf{v}}(W)=U_{2}$. Therefore, $\Sigma(\mathcal{B})$ is topologically mixing.

(ii)$\Rightarrow$(i). From Theorem \ref{theorem:3.4-0}, $\mathbb{H}_{n}$ and $\mathbb{V}_{n}$ are non-compressible for all $n\geq 2$. Then, for $n\geq 2$, any pattern in $\Sigma_{1\times n}(p)$ or $\Sigma_{n\times 1}(p)$ can be extended to $\mathbb{Z}^{2}$ by using the local patterns in $\mathcal{B}$. It can be easily verified that (ii)$\Rightarrow$(i); the details are omitted. The proof is complete. \hspace{0.5cm} $\square$
\medbreak

When $\mathbb{H}_{2}(\mathcal{B})$ or $\mathbb{V}_{2}(\mathcal{B})$ is degenerated, Theorems \ref{theorem:3.6} and \ref{theorem:3.7} can be generalized to the case when $\mathcal{B}$ satisfies locally corner-extendable conditions and local crisscross-extendability.

In the following, local crisscross-extendability is useful for extending a rectangular block to a global pattern.

When $\mathcal{B}\subset\Sigma_{2\times 2}(p)$ is not locally crisscross-extendable, $\mathcal{B}$ can be reduced to $\mathcal{B}_{c}\subseteq \mathcal{B}$ such that $\mathcal{B}_{c}$ is locally crisscross-extendable and $\Sigma(\mathcal{B}_{c})=\Sigma(\mathcal{B})$. The details are omitted here.

The following theorem shows that the local corner-extendable conditions and local crisscross-extendability imply rectangular extendability.

\begin{theorem}
\label{theorem:3.13}
If $\mathcal{B}$ satisfies either $C(1)$ and $C(3)$ or $C(2)$ and $C(4)$, then the following statements are equivalent.
\item[(i)] $\mathcal{B}$ is rectangle-extendable.

\item[(ii)]$\mathcal{B}$ is locally crisscross-extendable.
\end{theorem}

\textit{Proof.}
Clearly, (i) implies (ii).

(ii)$\Rightarrow$(i). Assume that $\mathcal{B}$ satisfies $C(1)$ and $C(3)$. The case in which it satisfies $C(2)$ and $C(4)$ is similar.
Let $U_{m\times n}\in\Sigma_{m \times n}(\mathcal{B})$, $m,n\geq 2$. Since $\mathcal{B}$ satisfies $C(1)$ and $C(3)$, from (ii), $U_{m\times n}$ can be extended in both positive and negative vertical directions by using the local patterns in $\mathcal{B}$, as follows.

\begin{equation*}
\psfrag{a}{{\footnotesize$n$}}
\psfrag{b}{{\footnotesize$m$}}
\psfrag{c}{{\footnotesize$n+2$}}
\psfrag{e}{.}
\includegraphics[scale=0.65]{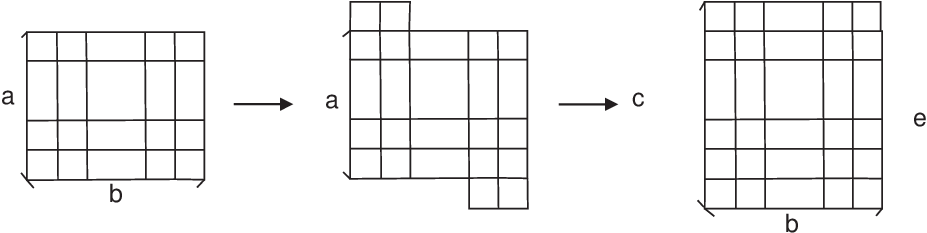}
\end{equation*}
\begin{equation*}
\text{Figure 3.6.}
\end{equation*}
Similarly, the above pattern can be extended in both positive horizontal and negative horizontal directions using the local patterns in $\mathcal{B}$. Therefore, by the above method, $U_{m\times n}$ can be extended to $\mathbb{Z}^{2}$ using the local patterns in $\mathcal{B}$. The proof is complete. \hspace{0.5cm} $\square$
\medbreak

Theorem \ref{theorem:3.7} can now be proven and generalized with a slight modification.
\begin{theorem}
\label{theorem:3.14}
If
\begin{enumerate}
\item[(i)] $\mathcal{B}\subset \Sigma_{2\times 2}(p)$ is locally crisscross-extendable, and

\item[(ii)] $\mathcal{B}$ satisfies three of the locally corner-extendable conditions $C(i)$, $1\leq i\leq 4$,
\end{enumerate}
then
$\mathbb{H}_{n}(\mathcal{B})$ and $\mathbb{V}_{n}(\mathcal{B})$ are weakly primitive for all $n\geq 2$ if and only if
$\Sigma(\mathcal{B})$ is topologically mixing.
\end{theorem}

\textit{Proof.}
$(\Rightarrow)$. From (ii), without loss of generality, assume that $\mathcal{B}$ satisfies conditions $C(1)$, $C(2)$ and $C(3)$.

Let $R_{1}$ and $R_{2}$ be finite sublattices of $\mathbb{Z}^{2}$. Since $\mathcal{B}$ satisfies $C(1)$ and $C(2)$, as in the proof of Theorem \ref{theorem:3.7}, there exists $M(R_{1},R_{2})\geq 1$ such that for all $\mathbf{v}=(v_{1},v_{2})\in\mathbb{Z}^{2}$ with $d(R_{1},R_{2}+\mathbf{v})\geq M$ and any two allowable patterns $U_{1}\in\Pi_{R_{1}}(\Sigma(\mathcal{B}))$ and $U_{2}\in\Pi_{R_{2}+\mathbf{v}}(\Sigma(\mathcal{B}))$, $U_{1}$ and $U_{2}$ can be extended as $U_{r}$ on the rectangular lattice using the local patterns in $\mathcal{B}$.

Since $\mathcal{B}$ is locally crisscross-extendable and satisfies conditions $C(1)$ and $C(3)$, by Theorem \ref{theorem:3.13}, $\mathcal{B}$ is rectangle-extendable.
Then, $U_{r}$ can be extended as $W\in\Sigma(\mathcal{B})$ with $\Pi_{R_{1}}(W)=U_{1}$ and $\Pi_{R_{2}+\mathbf{v}}(W)=U_{2}$. Therefore, $\Sigma(\mathcal{B})$ is topologically mixing.

$(\Leftarrow)$. From (i) and (ii), by Theorem \ref{theorem:3.13}, $\mathcal{B}$ is rectangle-extendable. Then, for $n\geq 2$, any pattern in $\Sigma_{2\times n}(\mathcal{B})$ can be extended to $\mathbb{Z}^{2}$ using the local patterns in $\mathcal{B}$. Therefore, the fact that $\Sigma(\mathcal{B})$ is topologically mixing implies that
$\mathbb{H}_{n}(\mathcal{B})$ is weakly primitive for all $n\geq 2$. Similarly, $\mathbb{V}_{n}(\mathcal{B})$ is weakly primitive for all $n\geq 2$.

The proof is complete. \hspace{0.5cm} $\square$
\medbreak

The following example demonstrates that the locally corner-extendable conditions in Theorem \ref{theorem:3.14} are crucial: if locally corner-extendable conditions are not satisfied, then local crisscross-extendability (or rectangle-extendability) and primitivity may not imply topological mixing.

\begin{example}
\label{example:3.15}

Let

\begin{equation*}
\mathcal{B}_{\pi/4}=\left\{
\begin{array}{c}
\psfrag{a}{$u_{2}$}
\psfrag{b}{$u_{3}$}
\psfrag{c}{$u_{1}$}
\psfrag{e}{$u_{4}$}
\psfrag{d}{$: u_{4}\geq u_{1}\text{ and } u_{1},u_{2},u_{3},u_{4}\in\{0,1\}$ }
\includegraphics[scale=0.8]{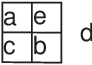}
\end{array}
\hspace{6.0cm}\right\},
\end{equation*}
which requires that diagonal lines with slope 1 are non-decreasing from left to right.

Clearly,

\begin{equation*}
\mathbb{H}_{2}(\mathcal{B}_{\pi/4})=\mathbb{V}_{2}(\mathcal{B}_{\pi/4})=
\left[
\begin{array}{cccc}
1 & 1 & 1 & 1 \\
1 & 1 & 1 & 1 \\
0& 1 & 0 & 1 \\
0 & 1 & 0 & 1
\end{array}
\right].
\end{equation*}
From (\ref{eqn:2.9}) and (\ref{eqn:2.10}), $\mathbb{H}_{n}$ is non-compressible for all $n\geq 2$. Since $\mathbb{V}_{2}=\mathbb{H}_{2}$, $\mathbb{V}_{n}$ is also non-compressible for all $n\geq 2$. By Theorem \ref{theorem:3.5}, $\mathcal{B}_{\pi/4}$ is rectangle-extendable. In particular, $\mathcal{B}_{\pi/4}$ is locally crisscross-extendable.

From the rule of $\mathcal{B}_{\pi/4}$, it can be easily proven that $\mathbb{H}_{n}^{n}=\mathbb{V}_{n}^{n}>0$ for all $n\geq 2$.
However, let $R_{1}=R_{2}=\mathbb{Z}_{2\times 2}$. Consider $U_{0}=\{0\}^{\mathbb{Z}^{2}}$ and $U_{1}=\{1\}^{\mathbb{Z}^{2}}$. Clearly, $U_{0},U_{1}\in\Sigma(\mathcal{B}_{\pi/4})$, but $\Pi_{R_{1}}(U_{1})$ cannot connect with $\Pi_{R_{2}+(i,i)}(U_{0})=\Pi_{\mathbb{Z}_{2\times 2}((i,i))}(U_{0})$ using the local patterns in $\mathcal{B}_{\pi/4}$ for all $i\geq 2 $. Then, $\Sigma(\mathcal{B}_{\pi/4})$ is not topologically mixing. Therefore, local crisscross-extendability and primitivity do not imply topological mixing. This claim does not contradict Theorem \ref{theorem:3.14} since $\mathcal{B}_{\pi/4}$ does not satisfy conditions $C(2)$ and $C(4)$:
neither
\begin{equation*}
\begin{array}{lcr}
 \psfrag{a}{{\small$0$}}
 \psfrag{b}{{\small$1$}}
  \psfrag{c}{$\in\Sigma_{\mathbb{L}_{2}}(\mathcal{B}_{\pi/4})$}
\includegraphics[scale=0.7]{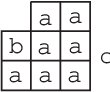} & \hspace{2.0cm}  \text{nor} &  \hspace{1.0cm}
 \psfrag{a}{{\small$0$}}
 \psfrag{b}{{\small$1$}}
  \psfrag{c}{$\in\Sigma_{\mathbb{L}_{4}}(\mathcal{B}_{\pi/4})$}
\includegraphics[scale=0.7]{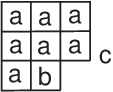}
\end{array}
\end{equation*}
can be extended to $\mathbb{Z}_{3\times 3}$ using the local patterns in $\mathcal{B}_{\pi/4}$.

\end{example}

\begin{remark}
\label{remark:3.5}
The non-degeneracy of $\mathbb{H}_{2}$ and $\mathbb{V}_{2}$ and the locally corner-extendable conditions are used to extend a single local pattern to be a global pattern, which are intrinsically different from mixing properties. Mixing is associated with two given local patterns that are parts of two global patterns. In fact, these extendability conditions alone cannot imply any mixing property. For example,

\begin{equation*}
\mathbb{H}_{2}(\mathcal{B})=\mathbb{V}_{2}(\mathcal{B})=\left[\begin{array}{cccc}  1 & 0 & 0 & 1 \\ 0 & 1 & 1 & 0 \\ 0& 1 & 1 & 0 \\ 1 & 0 & 0 & 1 \end{array}\right]
\end{equation*}
are non-degenerated, and $\Sigma(\mathcal{B})$ is not topologically mixing. Actually, $\mathbb{H}_{2}(\mathcal{B})$ and $\mathbb{V}_{2}(\mathcal{B})$ are not primitive.
\end{remark}

\numberwithin{equation}{section}

\section{Invariant diagonal cycles and primitive commutative cycles}
\label{sec:4}
\hspace{0.5cm}
According to Propositions \ref{proposition:2.2} and \ref{proposition:2.3}, the recursive formula from a higher-order transition matrix to a lower-order transition matrix using the connecting operator can be used to introduce invariant diagonal cycles and primitive commutative cycles that enable the connecting operator to be used to provide finitely sufficient conditions for the primitivity of $\mathbb{H}_{n}$ and $\mathbb{V}_{n}$ for $n\geq 2$. For brevity, only $\mathbb{H}_{n}$ is considered here. The discussion for $\mathbb{V}_{n}$ is similar to that for $\mathbb{H}_{n}$.

\subsection{Invariant diagonal cycles}
\label{sec:4-1}
This subsection introduces invariant diagonal cycles of the connecting operator to provide finitely sufficient conditions for the primitivity of $\mathbb{H}_{n}$ or $\mathbb{V}_{n}$.

First, the diagonal index set is defined by
\begin{equation*}
\mathcal{D}_{p}=\left\{1+j(p+1)| j\in\mathcal{S}_{p}\right\}.
\end{equation*}
Clearly, if $\beta_{1},\beta_{2},\cdots ,\beta_{q+1}\in\mathcal{D}_{p}$, then $H_{m,n+q;\beta_{1};\beta_{2};\cdots ;\beta_{q+1}}$ lies on the diagonal of $\mathbb{H}_{n+q}^{m}$ in (\ref{eqn:2.34}).

\begin{definition}
\label{definition:4.9}
\item[(i)] For $q\geq 1$, a finite sequence  $\overline{\beta}_{q}=\beta_{1}\beta_{2}\cdots\beta_{q}\beta_{1}$ is called a diagonal cycle with length $q$ if $\beta_{j}\in \mathcal{D}_{p}$ for $1\leq j\leq q$.

\item[(ii)] A diagonal cycle $\overline{\beta}_{q}=\beta_{1}\beta_{2}\cdots\beta_{q}\beta_{1}$ is called an $S$-invariant diagonal cycle of order $(m,q)$ if there exist $m\geq 2$ and an invariant index set  $\mathcal{K}\subseteq \left\{1,2,\cdots,p^{m-1}\right\}$ such that

\begin{equation}\label{eqn:4.9}
\underset{k\in\mathcal{K}}{\sum}\left(S_{m;\beta_{1},\beta_{2}}S_{m;\beta_{2},\beta_{3}}\cdots S_{m;\beta_{q},\beta_{1}}\right)_{k,l}\geq 1
\end{equation}
for all $l\in\mathcal{K}$.

\item[(iii)] A $W$-invariant diagonal cycle can be defined analogously.
\end{definition}

Notably, it can be easily shown that for any $n\geq 1$,

\begin{equation}\label{eqn:4.10}
\underset{k\in\mathcal{K}}{\sum}\left((S_{m;\beta_{1},\beta_{2}}S_{m;\beta_{2},\beta_{3}}\cdots S_{m;\beta_{q},\beta_{1}})^{n}\right)_{k,l}\geq 1
\end{equation}
for all $l\in\mathcal{K}$ if (\ref{eqn:4.9}) holds. The case for $W$-invariant diagonal cycles is similar.

The following notation is used in proving the theorem for the primitivity of $\mathbb{H}_{n}$.

\begin{definition}
\label{definition:4.8}
Let $\mathbb{M}=\left[M_{i,j}\right]_{N\times N}$, where $M_{i,j}$ is an $M\times M$ non-negative matrix for $1\leq i,j\leq N$. The indicator matrix $\Lambda(\mathbb{M})=[m_{i,j}]_{N\times N}$ of $\mathbb{M}$ is defined by

\begin{equation*}
\left\{
\begin{array}{ll}
m_{i,j}=1 & \text{if }   |M_{i,j}|>0 ,\\
m_{i,j}=0 & \text{otherwise, }
\end{array}
\right.
\end{equation*}
where $|M_{i,j}|$ is the sum of all entries in $M_{i,j}$.
\end{definition}

The following lemma is essential for establishing the primitivity of $\mathbb{H}_{n}$ using the invariant diagonal cycle.

\begin{lemma}
\label{lemma:4.9}
Suppose $\mathbb{M}=\left[M_{i,j}\right]_{N\times N}$, where $M_{i,j}$ is an $M\times M$ non-negative matrix for $1\leq i,j\leq N$. Let $\Lambda(\mathbb{M})=[m_{i,j}]_{N\times N}$ be the indicator matrix of $\mathbb{M}$.
Let $\mathbb{M}^{n}=\left[M_{n;i,j}\right]_{N\times N}$ for $n\geq 1$. If
\begin{enumerate}
\item[(i)] $\Lambda(\mathbb{M})$ is primitive,
\item[(ii)] $M_{i,j}$ is either non-compressible or zero, $1\leq i,j\leq N$, and
\item[(iii)] there exist $n\geq 1$ and $1\leq k \leq N$ such that $M_{n;k,k}$ is primitive,
\end{enumerate}
then $\mathbb{M}$ is primitive.
\end{lemma}

\textit{Proof.}
Since $\Lambda(\mathbb{M})$ and $M_{n;k,k}$ are primitive, there exists $N_{1}\geq 1$ such that $\Lambda(\mathbb{M})^{N_{1}}>0$ and $M_{n;k,k}^{N_{1}}>0$. By (ii), for any $l\geq 1$ and $1\leq i,j\leq N$, if $\left(\Lambda(\mathbb{M})^{l}\right)_{i,j}>0$, then $M_{l;i,j}$ is non-compressible.

Take $N_{2}=(2+n) N_{1}$. The fact that for any $1\leq i,j\leq N$,

\begin{equation*}
M_{N_{2};i,j}\geq M_{N_{1};i,k}M_{n;k,k}^{N_{1}}M_{N_{1};k,j}>0,
\end{equation*}
can be easily seen.
Therefore, $\mathbb{M}$ is primitive. \hspace{0.5cm} $\square$
\medbreak

When an invariant diagonal cycle of order $(m,q)$ exists, the following theorem shows that the primitivity of $\underset{l\in\mathcal{K}}{\sum}H_{m,n;\beta_{1}}^{(l)}$ up to finite order $q+1$ implies the primitivity of $\mathbb{H}_{n}$ for all $n\geq 2$.

\begin{theorem}
\label{theorem:4.11-1}
Given $\mathcal{B}\subset\Sigma_{2 \times 2}(p)$, if
\begin{enumerate}
\item[(i)] $\mathbb{H}_{2}(\mathcal{B})$ is non-degenerated,

\item[(ii)] there exists an $S$-invariant diagonal cycle $\overline{\beta}_{q}=\beta_{1}\beta_{2}\cdots\beta_{q}\beta_{1}$ of order $(m,q)$ with its invariant index set $\mathcal{K}$, and

\item[(iii)] $\underset{l\in\mathcal{K}}{\sum}H_{m,n;\beta_{1}}^{(l)}$ is primitive for $2\leq n\leq q+1$,
\end{enumerate}
then $\mathbb{H}_{n}$ is primitive for all $n\geq 2$.
\end{theorem}

\textit{Proof.}
The result that $\mathbb{H}_{n}$ is primitive for $n\geq 2 $ is proven by induction, as follows. For $s\geq 0$,
the statement $P(s)$ means that $\mathbb{H}_{n}$ is primitive for $sq+2 \leq n\leq (s+1)q+1$.

For $2\leq n\leq q+1$, let $\mathbb{H}_{n}=[H_{n;\alpha}]_{p\times p}$. By Theorem \ref{theorem:3.4-0}, $H_{n;\alpha}$ is non-compressible for $1\leq \alpha\leq p^{2}$. Clearly, the indicator matrix of $\mathbb{H}_{n}=[H_{n;\alpha}]_{p\times p}$ is $p\times p$ full matrix. From (iii),

\begin{equation*}
H_{m,n;\beta_{1}}= \underset{l=1}{\overset{p^{m-1}}{\sum}}H_{m,n;\beta_{1}}^{(l)}\geq \underset{l\in\mathcal{K}}{\sum}H_{m,n;\beta_{1}}^{(l)}
 \end{equation*}
is primitive and is on the diagonal of $\mathbb{H}_{n}^{m}=[H_{m,n;\alpha}]_{p\times p}$. Hence, by Lemma \ref{lemma:4.9}, $\mathbb{H}_{n}$ is primitive for $2\leq n\leq q+1$ and $P(0)$ is true..

Assume that $P(t)$ follows for some $t\geq 0$, meaning that,  $\mathbb{H}_{n}$ is primitive for $tq+2 \leq n\leq (t+1)q+1$.

For $(t+1)q+2 \leq n\leq (t+2)q+1$, let $n=(t+1)q+r$, where $2\leq r\leq q+1$. Let $N=(t+1)q+1$, define $\overline{\beta}_{N-1}=\left(\beta_{1}\beta_{2}\cdots\beta_{q}\right)^{t+1}\beta_{1}$. From (\ref{eqn:4.10}), $\overline{\beta}_{N-1}$ is an $S$-invariant diagonal cycle of order $(m,N-1)$ with invariant index set $\mathcal{K}$.

From (\ref{eqn:2.34}), let $\mathbb{H}_{n}=[H_{n;i,j}]_{p^{N}\times p^{N}}$. Then, from (\ref{eqn:2.19}),

\begin{equation*}
\mathbb{H}_{n}=[H_{n;i,j}]_{p^{N}\times p^{N}}=\left( \mathbb{H}_{N} \right)_{p^{N}\times p^{N}}\circ\left[
E_{p^{N-1}\times p^{N-1}}\otimes \left[H_{r;i,j}\right]_{p\times p}
\right].
\end{equation*}
By Theorem \ref{theorem:3.4-0}, $H_{r;i,j}$ is non-compressible, $1\leq i,j\leq p$. Then, $\mathbb{H}_{N}$ is the indicative matrix of $\mathbb{H}_{n}=[H_{n;i,j}]_{p^{N}\times p^{N}}$ and $H_{n;i,j}$ is either non-compressible or zero for $1\leq i,j\leq p^{N}$. By the assumption for $P(t)$, $\mathbb{H}_{N}$ is primitive.

Let $\mathbb{H}^{m}_{n}=[H_{m,n;i,j}]_{p^{N}\times p^{N}}=[H_{m,n;\alpha_{1};\alpha_{2};\cdots;\alpha_{N}}]_{p^{N}\times p^{N}}$.
From (\ref{eqn:2.36}),

\begin{equation*}
\begin{array}{rl}
H_{m,n;\overline{\beta}_{N-1}} \equiv & H_{m,n;\beta_{1};\beta_{2};\cdots ;\beta_{\bar{q}};\cdots;\beta_{1};\beta_{2};\cdots ;\beta_{q};\beta_{1}} \\
 & \hspace{1.0cm}\underset{(t+1)\text{ times}}{ \underbrace{     \hspace{3.2cm}              } }\\
 & \\
= &
\underset{k,l=1}{\overset{p^{m-1}}{\sum}}
((S_{m;\beta_{1},\beta_{2}}S_{m;\beta_{2},\beta_{3}}\cdots S_{m;\beta_{q},\beta_{1}})^{t+1})_{k,l}
H_{m,r;\beta_{1}}^{(l)}\\
\geq & \underset{l\in\mathcal{K}}{{\sum}}H_{m,r;\beta_{1}}^{(l)}.
\end{array}
\end{equation*}
$H_{m,n;\overline{\beta}_{N-1}}$ is on the diagonal of $\mathbb{H}_{n}^{m}$. By Lemma \ref{lemma:4.9},
$\mathbb{H}_{n}$ is primitive for  $(t+1)\bar{q}+2 \leq n\leq (t+2)\bar{q}+1$, so $P(t+1)$ holds.

Therefore, by induction, $P(s)$ is true for all $s\geq 0$, implying that $\mathbb{H}_{n}$ is primitive for all $n\geq 2$. The proof is complete. \hspace{0.5cm} $\square$
\medbreak

The following example illustrates the application of Theorem \ref{theorem:4.11-1}.

\begin{example}
\label{example:4.12}
Consider
\begin{equation*}
\mathbb{H}_{2}(\mathcal{B})=
\left[
\begin{array}{cccc}
1 & 0 & 0 & 1 \\
1 & 1 & 1 & 0 \\
1& 0 & 0 & 1 \\
0 & 1 & 1 & 0
\end{array}
\right].
\end{equation*}
Clearly, $\mathbb{H}_{2}$ is non-degenerated. From (\ref{eqn:2.27}),
\begin{equation*}
S_{3;1,1}=C_{3;1,1}=
\left[
\begin{array}{cccc}
1 & 0  & 0 & 0 \\
0 & 0  & 0 & 0 \\
0 & 0  & 0 & 1 \\
0 & 0  & 1 & 0
\end{array}
\right] .
\end{equation*}
Let $\overline{\beta}_{1}=11$ and $\mathcal{K}=\{3,4\}$. Since
\begin{equation*}
\underset{k\in\mathcal{K}}{\sum}\left(S_{3;1,1}\right)_{k,l}\geq 1
\end{equation*}
for $l\in\mathcal{K}$,
$\overline{\beta}_{1}$ is an S-invariant diagonal cycle of order $(3,1)$ with index set $\mathcal{K}$.
Clearly,

\begin{equation*}
\underset{l\in\mathcal{K}}{\sum}H_{3,2;1}^{(l)} =H_{2;1,2}H_{2;2,1}H_{2;1,1}+H_{2;1,2}H_{2;2,2}H_{2;2,1}
     =   \left[
\begin{array}{cc}
2 & 1  \\
1 & 1
\end{array}
\right]
\end{equation*}
is primitive. From Theorem \ref{theorem:4.11-1}, $\mathbb{H}_{n}$ is primitive for all $n\geq 2$.

The fact that $\mathbb{V}_{2}(\mathcal{B})$ does not have invariant diagonal cycle up to $m=7$ can be verified, but whether $\mathbb{V}_{2}(\mathcal{B})$ has an invariant diagonal cycle for larger $m$ is unclear. To deal with this difficulty, the following subsection introduces another criterion for establishing primitivity using primitive commutative cycles.
The topological mixing of $\Sigma(\mathcal{B})$ is proven in Example \ref{example:5.8}.

\end{example}


\subsection{Primitive commutative cycles}
\label{sec:4-2}
\setcounter{equation}{10}
This subsection introduces primitive commutative cycles to obtain another finitely  sufficient condition for the primitivity of $\mathbb{H}_{n}$ or $\mathbb{V}_{n}$ when invariant diagonal cycles are not available.

For $q,q'\geq 1$, let $I_{q}=i_{1}i_{2}\cdots i_{q}i_{1}$ and $J_{q'}=j_{1}j_{2}\cdots j_{q'}j_{1}$ be two cycles, where $i_{k},j_{l}\in\{1,2,\cdots,p\}$ for $1\leq k\leq q$ and $1\leq l\leq q'$.

\begin{definition}
\label{definition:5.1}
If $j_{1}=i_{1}$, let $(I_{q}J_{q'})=i_{1}i_{2}\cdots i_{q}i_{1}j_{2}\cdots j_{q'}i_{1}$ and $(J_{q'}I_{q})=i_{1}j_{2}\cdots j_{q'}i_{1}i_{2}\cdots i_{q}i_{1}$. The pair $(I_{q}J_{q'})$ and $(J_{q'}I_{q})$ is called  a commutative cycle pair.
\end{definition}

Given a commutative cycle pair $(I_{q}J_{q'})$ and $(J_{q'}I_{q})$, denote the index of $(I_{q}J_{q'})$ and $(J_{q'}I_{q})$ by $\langle m,\bar{\alpha};K,L\rangle$, where

\begin{equation}\label{eqn:5.1}
\left\{
\begin{array}{l}
m=q+q' \\
\bar{\alpha}=\psi(i_{1}-1,i_{1}-1)\\
K=\psi(i_{2}-1,\cdots ,i_{q}-1,i_{1}-1,j_{2}-1,\cdots ,j_{q'}-1)\\
L=\psi(j_{2}-1,\cdots ,j_{q'}-1,i_{1}-1,i_{2}-1,\cdots ,i_{q}-1).
\end{array}
\right.
\end{equation}

From (\ref{eqn:2.20}), it is easy to check that

\begin{equation}\label{eqn:5.2}
\left\{
\begin{array}{l}
H_{n;i_{1},i_{2}}H_{n;i_{2},i_{3}}\cdots H_{n;i_{q},i_{1}}H_{n;i_{1},j_{2}}H_{n;j_{2},j_{3}}\cdots H_{n;j_{q'},i_{1}}=H_{m,n;\bar{\alpha}}^{(K)} \\
\\
H_{n;i_{1},j_{2}}H_{n;j_{2},j_{3}}\cdots H_{n;j_{q'},i_{1}}H_{n;i_{1},i_{2}}H_{n;i_{2},i_{3}}\cdots H_{n;i_{q},i_{1}}=H_{m,n;\bar{\alpha}}^{(L)}.
\end{array}
\right.
\end{equation}
The number $\bar{\alpha}$ is a member of the diagonal index set $\mathcal{D}_{p}$, and then $H_{m,n;\bar{\alpha}}$ lies on the diagonal of $\mathbb{H}_{n}^{m}$.

\begin{definition}
\label{definition:5.2}
 A commutative cycle pair $(I_{q}J_{q'})$ and $(J_{q'}I_{q})$ with index $\langle m,\bar{\alpha};K,L\rangle$ is called an $H$-primitive commutative cycle pair if
$H_{m,2;\bar{\alpha}}^{(K)}$ and $H_{m,2;\bar{\alpha}}^{(L)}$ are primitive.

\end{definition}
A $V$-primitive commutative cycle pair is similarly specified, and the details are omitted. The commutative cycle pair can compensate for each other, and they can establish the primitivity of $\mathbb{H}_{n}$ for all $n\geq 2$.
Indeed, the following theorem provides a sufficient condition for the primitivity of $\mathbb{H}_{n}$ when $\mathbb{H}_{2}$ is non-degenerated. Similar results hold for $\mathbb{V}_{n}$.

 \begin{theorem}
\label{theorem:5.3}
Given $\mathcal{B}\subset\Sigma_{2\times 2}(p)$, if
\begin{enumerate}
\item[(i)] $\mathbb{H}_{2}$ is non-degenerated, and

\item[(ii)] there exists an $H$-primitive commutative cycle pair $(I_{q}J_{q'})$ and $(J_{q'}I_{q})$ with index $\langle m,\bar{\alpha};K,L\rangle$ such that $(S_{m;\bar{\alpha},\bar{\alpha}})_{K,L}=1$ or $(S_{m;\bar{\alpha},\bar{\alpha}})_{L,K}=1$,
\end{enumerate}
then $\mathbb{H}_{n}$ is primitive for all $n\geq 2$.
\end{theorem}

\textit{Proof.}
Suppose $(S_{m;\bar{\alpha},\bar{\alpha}})_{K,L}=1$. The case for $(S_{m;\bar{\alpha},\bar{\alpha}})_{L,K}=1$ is similar.

First, $H_{m,n;\bar{\alpha}}^{(K)}$ and $H_{m,n;\bar{\alpha}}^{(L)}$ will be shown by induction to be primitive for $n\geq 2$ by induction.
From (ii), the primitivity of $H_{m,2;\bar{\alpha}}^{(K)}$ and $H_{m,2;\bar{\alpha}}^{(L)}$ holds for $n=2$. Assume that this result holds for $n=t$, $t\geq 2$.

 Let $H_{m,2;\bar{\alpha}}^{(K)}=\left[H_{m,2;\bar{\alpha};\alpha}^{(K)}\right]_{p\times p}$. From (\ref{eqn:2.33}),

\begin{equation}\label{eqn:4.100}
H_{m,2;\bar{\alpha};\alpha}^{(K)}=\underset{l=1}{\overset{p^{m-1}}{\sum}}(S_{m;\bar{\alpha},\alpha})_{K,l}
\end{equation}
for all $1\leq i,j\leq p$. Let $\Lambda=\Lambda\left(H_{m,2;\bar{\alpha}}^{(K)}\right)$ be the indicator matrix of $H_{m,2;\bar{\alpha}}^{(K)}=\left[H_{m,2;\bar{\alpha};\alpha}^{(K)}\right]_{p\times p}$.
The primitivity of $H_{m,2;\bar{\alpha}}^{(K)}$ implies $\Lambda$ is primitive.

Consider the case for $n=t+1$. Let $H_{m,t+1;\bar{\alpha}}^{(K)}=\left[H_{m,t+1;\bar{\alpha};\alpha}^{(K)}\right]_{p\times p}$, by Proposition \ref{proposition:2.2},

\begin{equation}\label{eqn:4.101}
H_{m,t+1;\bar{\alpha};\alpha}^{(K)}=\underset{l=1}{\overset{p^{m-1}}{\sum}}(S_{m;\bar{\alpha},\alpha})_{K,l}
H_{m,t;\alpha}^{(l)}
\end{equation}
for all $1\leq \alpha\leq p^{2}$. By Theorem \ref{theorem:3.6}, every pattern $U_{m\times 2}\in\Sigma_{m\times 2}(\mathcal{B})$ can be extended to $\mathbb{Z}_{m\times 3}$ by using the local patterns in $\mathcal{B}$. Thus, if $(S_{m;\bar{\alpha},\alpha})_{K,l}=1$, then $H_{m,t;\alpha}^{(l)}$ is not a zero matrix for $1\leq \alpha\leq p^{2}$ and $1\leq l\leq p^{m-1}$. Hence,
$\Lambda\left(H_{m,2;\bar{\alpha}}^{(K)}\right)$ is also the indicator matrix of $H_{m,t+1;\bar{\alpha}}^{(K)}=\left[H_{m,t+1;\bar{\alpha};\alpha}^{(K)}\right]_{p\times p}$. Moreover, from (\ref{eqn:4.101}) and Theorem \ref{theorem:3.4-0}, $H^{(K)}_{m,t+1;\bar{\alpha};\alpha}$ is either non-compressible or zero, $1\leq \alpha\leq p^{2}$. From (\ref{eqn:4.101}), $H^{(K)}_{m,t+1;\bar{\alpha};\bar{\alpha}}$ is primitive and on the diagonal of $H_{m,t+1;\bar{\alpha}}^{(K)}$. From Lemma \ref{lemma:4.9}, $H_{m,t+1;\bar{\alpha}}^{(K)}$ is primitive.

Let $A$ and $B$ be non-negative and non-compressible matrices. If $AB$ is primitive, then $BA$ can be easily verified also to be primitive. By (\ref{eqn:5.2}), $H_{m,t+1;\bar{\alpha}}^{(L)}$ is primitive.
Hence, the case for $n=t+1$ holds. Therefore, $H_{m,n;\bar{\alpha}}^{(K)}$ and $H_{m,n;\bar{\alpha}}^{(L)}$ are primitive for $n\geq 2$.

Now, $\mathbb{H}_{n}$ will be shown to be primitive for all $n\geq 2$. In the case $n=2$, let $\mathbb{H}_{2}=[H_{2;\alpha}]_{p\times p}$. By Theorem \ref{theorem:3.4-0}, $H_{2;\alpha}$ is non-compressible for $1\leq \alpha\leq p^{2}$. The indicator matrix of $\mathbb{H}_{2}=[H_{2;\alpha}]_{p\times p}$ is a $p\times p$ full matrix. From (ii),

\begin{equation*}
H_{m,2;\bar{\alpha}}= \underset{l=1}{\overset{p^{m-1}}{\sum}}H_{m,2;\bar{\alpha}}^{(l)}\geq H_{m,2;\bar{\alpha}}^{(K)}
 \end{equation*}
is primitive and is on the diagonal of $\mathbb{H}_{2}^{m}=[H_{m,2;\alpha}]_{p\times p}$. Then, by Lemma \ref{lemma:4.9}, $\mathbb{H}_{2}$ is primitive.

For $n\geq 3$, from (\ref{eqn:2.19}),

\begin{equation*}
\mathbb{H}_{n}=\left[H_{n;i,j}\right]_{p^{2}\times p^{2} }=\left( \mathbb{H}_{2} \right)_{p^{2}\times p^{2}}\circ\left[
E_{p\times p}\otimes \left[H_{n-1;\alpha}\right]_{p\times p}
\right].
\end{equation*}

From (ii), by Theorem \ref{theorem:3.4-0}, if $\left( \mathbb{H}_{2} \right)_{i,j}=1$, then $H_{n;i,j}$ is not a zero matrix. Hence, $\mathbb{H}_{2}$ is the indicator matrix of $\mathbb{H}_{n}=\left[H_{n;i,j}\right]_{p^{2}\times p^{2} }$.

Let $\mathbb{H}_{n}^{m}=\left[H_{m,n;\alpha_{1};\alpha_{2}}\right]_{p^{2}\times p^{2} }$.
Since $(S_{m;\bar{\alpha},\bar{\alpha}})_{K,L}=1$,

\begin{equation*}
H_{m,n;\bar{\alpha};\bar{\alpha}}
=
\underset{k,l=1}{\overset{p^{m-1}}{\sum}}
(S_{m;\bar{\alpha};\bar{\alpha}})_{k,l}
H_{m,n-1;\bar{\alpha}}^{(l)}
\geq H_{m,n-1;\bar{\alpha}}^{(L)}.
\end{equation*}
Since $H_{m,n-1;\bar{\alpha}}^{(L)}$ is primitive, $H_{m,n;\bar{\alpha};\bar{\alpha}}$ is primitive.
Notably, $H_{m,n;\bar{\alpha};\bar{\alpha}}$ is on the diagonal of $\mathbb{H}_{n}^{m}$. Therefore, from Lemma \ref{lemma:4.9}, $\mathbb{H}_{n}$ is primitive for all $n\geq 3$. The proof is complete. \hspace{0.5cm} $\square$
\medbreak

$\mathbb{H}$ (or $\mathbb{V}$) may have an invariant diagonal cycle and $\mathbb{V}$ (or $\mathbb{H}$) may have primitive commutative cycles. Therefore, combining these two conditions for the primitivity, the following theorem provides a finitely sufficient condition for topological mixing of $\Sigma(\mathcal{B})$ when $\mathbb{H}_{2}$ and $\mathbb{V}_{2}$ are non-degenerated.

\begin{theorem}
\label{theorem:5.7}
Given $\mathcal{B}\subset\Sigma_{2 \times 2}(p)$, if
\begin{enumerate}
\item[(i)] $\mathbb{H}_{2}(\mathcal{B})$ and $\mathbb{V}_{2}(\mathcal{B})$ are non-degenerated, and
\item[(ii)] $\mathcal{B}$ satisfies the conditions of Theorem \ref{theorem:4.11-1} or \ref{theorem:5.3} for $\mathbb{H}_{n}$ and $\mathbb{V}_{n}$,
\end{enumerate}
then $\Sigma(\mathcal{B})$ is topologically mixing.
\end{theorem}

\textit{Proof.}
Combining Theorems 3.6, 4.4 and 4.8 yields the result immediately.
\hspace{0.5cm} $\square$
\medbreak

Example 4.5, above, illustrates the application of Theorem \ref{theorem:5.7}.

\begin{example}
\label{example:5.8} (continued)

In Example \ref{example:4.12}, $\mathbb{H}_{2}(\mathcal{B})$ has an invariant diagonal cycle that satisfies condition (iii) of Theorem \ref{theorem:4.11-1}. However, $\mathbb{V}_{2}(\mathcal{B})$ does not have such an invariant diagonal cycle when $m\leq 7$, but it does have a primitive commutative cycle pair with $m=7$ that satisfies condition (ii) of Theorem \ref{theorem:5.3} as follows.
\begin{equation*}
\mathbb{V}_{2}(\mathcal{B})=
\left[
\begin{array}{cccc}
1 & 0 & 1 & 1 \\
0 & 1 & 1 & 0 \\
1& 0 & 0 & 1 \\
0 & 1 & 1 & 0
\end{array}
\right].
\end{equation*}
Clearly, $\mathbb{V}_{2}(\mathcal{B})$ is non-degenerated.

Let $I_{5}=211212$ and $J_{2}=222$. That

\begin{equation*}
V^{(12)}_{7,2;4}=V_{2;2,1}V_{2;1,1}V_{2;1,2}V_{2;2,1}V_{2;1,2}V_{2;2,2}V_{2;2,2}
=
\left[
\begin{array}{cc}
2 & 1 \\
1 & 1
\end{array}
\right]
\end{equation*}
and

\begin{equation*}
V^{(51)}_{7,2;4}=V_{2;2,2}V_{2;2,2}V_{2;2,1}V_{2;1,1}V_{2;1,2}V_{2;2,1}V_{2;1,2}
=
\left[
\begin{array}{cc}
2 & 1 \\
1 & 1
\end{array}
\right]
\end{equation*}
are primitive can be easily verified. Hence, $(I_{5}J_{2})$ and $(J_{2}I_{5})$ form a $V$-primitive commutative cycle pair with index $\langle 7,4;12,51\rangle$. Moreover,
\begin{equation*}
\left(W_{7;4,4}\right)_{12,51}=1.
\end{equation*}

Therefore, combining the result in Example \ref{example:4.12} and by Theorem \ref{theorem:5.7}, $\Sigma(\mathcal{B})$ is topologically mixing.
\end{example}

Under the local crisscross-extendibility and local corner-extendable conditions, Theorems \ref{theorem:4.11-1}, \ref{theorem:5.3} and \ref{theorem:5.7} can be generalized to some degenerated cases in which $\mathbb{H}_{2}$ or $\mathbb{V}_{2}$ contains zero rows or columns.
For completeness, the weakly non-degenerated case is introduced below.

\begin{definition}
\label{definition:4.53}
Given $\mathcal{B}\subset\Sigma_{2\times 2}(p)$, $\mathbb{H}_{2}(\mathcal{B})=[H_{2;i,j}]_{p\times p}$ is weakly non-degenerated if
\begin{enumerate}
\item[(i)] when both $H_{2;i,j_{1}}$ and $H_{2;i,j_{2}}$ are not zero matrices, $1\leq i,j_{1},j_{2}\leq p$,
\begin{equation*}
r(H_{2;i,j_{1}})=r(H_{2;i,j_{2}}), \text{and}
\end{equation*}

\item[(ii)] when both $H_{2;i_{1},j}$ and $H_{2;i_{2},j}$ are not zero matrices, $1\leq i_{1},i_{2},j \leq p$,
\begin{equation*}
c(H_{2;i_{1},j})=c(H_{2;i_{2},j}).
\end{equation*}
\end{enumerate}
Weak non-degeneracy of $\mathbb{V}_{2}(\mathcal{B})$ is defined analogously.
\end{definition}
Similar to Theorem \ref{theorem:3.6}, if $\mathcal{B}$ is locally crisscross-extendable and $\mathbb{H}_{2}$ and $\mathbb{H}_{2}$ are weakly non-degenerated, then $\mathcal{B}$ satisfies the local corner-filling conditions C(1), C(2) and C(4). For brevity, the details of the proof are omitted.

The following definition is introduced to enable the primitivity of compressible matrices to be easily expressed.

\begin{definition}
\label{definition:4.70}

If $A=[a_{i,j}]_{n\times n}$ is a matrix with $a_{i,j}\in\{0,1\}$, the associated saturated matrix $\mathbb{E}(A)=[e_{i,j}]_{n\times n}$ of $A$ is defined by

\begin{equation}\label{eqn:4.6}
\left\{
\begin{array}{rl}
e_{i,j}=0 & \hspace{1.0cm}\text{if }\underset{k=1}{\overset{n}{\sum}}a_{i,k}=0\text{ or }\underset{k=1}{\overset{n}{\sum}}a_{k,j}=0 , \\
& \\
e_{i,j}=1 &\hspace{1.0cm} \text{otherwise.} \
\end{array}
\right.
\end{equation}
\end{definition}
Clearly, given $A=[a_{i,j}]_{n\times n}$ with $a_{i,j}\in\{0,1\}$, if there exists $N\geq 1$ such that $A^{N}\geq \mathbb{E}(A)$, then $A$ is weakly primitive (weakly $N$-primitive); here, if $B=[b_{i,j}]_{n\times n}$ and $C=[c_{i,j}]_{n\times n}$ are two matrices, $B\geq C$ means $b_{i,j}\geq c_{i,j}$ for all $1\leq i,j\leq n$.

The following theorem can be obtained by an argument similar to that used in the non-degenerated case; the details of the proof are omitted. Notably, shift spaces such as the Golden Mean shift space can be applied to the following theorem.

\begin{theorem}
\label{theorem:4.210}
Given $\mathcal{B}\subset\Sigma_{2 \times 2}(p)$, if
\begin{enumerate}
\item[(i)] $\mathbb{H}_{2}(\mathcal{B})$ is weakly non-degenerated,

\item[(ii)] $\mathcal{B}$ is locally crisscross-extendable,

\item[(iii)] there exists an $S$-invariant diagonal cycle $\overline{\beta}_{q}=\beta_{1}\beta_{2}\cdots\beta_{q}\beta_{1}$ of order $(m,q)$ with its invariant index set $\mathcal{K}$,

\item[(iv)]  for $2\leq n\leq q+1$, there exists $a=a(n)\geq 1$ such that
\begin{equation*}
\left(\underset{l\in\mathcal{K}}{\sum}H_{m,n;\beta_{1}}^{(l)}\right)^{a}\geq  \mathbb{E}\left(H_{n;\beta_{1}}\right), \text{ and}
\end{equation*}
 \item[(v)] $\mathbb{H}_{n}$ is weakly primitive for $2\leq n\leq q+1$,

\end{enumerate}
then $\mathbb{H}_{n}$ is weakly primitive for all $n\geq 2$. Moreover, if $\mathbb{V}_{2}(\mathcal{B})$ also satisfies the conditions similar to (i)$\sim$(v), then $\Sigma(\mathcal{B})$ is topologically mixing.
\end{theorem}

\numberwithin{equation}{section}

\section{Strong specification}

\label{sec:6}
\hspace{0.5cm}

This section introduces the $k$ hole-filling condition ($($HFC$)_{k}$), and provides finitely sufficient conditions for the strong specification of $\Sigma(\mathcal{B})$.

 The main idea of finding sufficient conditions for strong specification is presented as follows. Clearly, strong specification is stronger than topological mixing. Apart from the processes in Fig. 3.1, which are associated with the situation in which regions $R_{1}$ and $R_{2}+\mathbf{v}$ are far away, the case in which one pattern is enclosed in another pattern, as in Fig. 5.1, must be studied.

\begin{equation*}
\psfrag{a}{{\footnotesize $U_{1}$}}
\psfrag{b}{{\footnotesize$U_{2}$}}
\psfrag{c}{}
\includegraphics[scale=1.0]{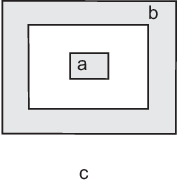}
\end{equation*}
\begin{equation*}
\text{Figure 5.1.}
\end{equation*}

Notably, in a study of topological mixing, Fig. 5.1 does not need to be considered because the separation distance can be chosen to be sufficiently large. However, in studying strong specification, $U_{1}$ and $U_{2}$ cannot be removed since the relative positions of $R_{1}$ and $R_{2}$ are fixed. Now, the sufficient condition is imposed to ensure that the gluing of $U_{1}$ and $U_{2}$ can be completed by the following two processes.

\begin{enumerate}
\item[Step (S-1):] Extend $U_{1}$ horizontally and vertically to form a crisscross pattern that touches $U_{2}$, as presented in Fig. 5.2.
\end{enumerate}

\begin{enumerate}
\item[Step (S-2):] Fill the holes that are surrounded by the rectangularly annular lattice to form a rectangular pattern, as presented in Fig. 5.3.
\end{enumerate}

\begin{equation*}
\begin{array} {cccccc}
\psfrag{a}{{\footnotesize  $U_{1}$}}
\psfrag{b}{{\footnotesize $U_{2}$}}
\psfrag{c}{{\tiny $(S-1)$}}
\includegraphics[scale=1.0]{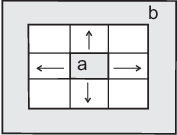}&
& & & &
\hspace{1.0cm}\psfrag{a}{{\tiny $(S-2)$}}
\includegraphics[scale=1.0]{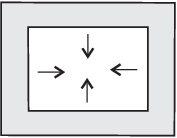}
\\&&&&&\\
\text{Figure 5.2.} && & &  & \hspace{1.0cm}\text{Figure 5.3.}
\end{array}
\end{equation*}
Then, repeat Step $(3)$ for topological mixing in Section 4 to extend the rectangular pattern to a global pattern on $\mathbb{Z}^{2}$.

The hole-filling condition in Step (S-2) is closely related to the extension property called square filling \cite{40-1,40-2}. In the following, the hole-filling condition is introduced.

First, for $M,N\geq 1$ and $i,j\in\mathbb{Z}$, the rectangularly annular lattice $\mathcal{A}_{M \times N;d}((i,j))$ with hole $\mathbb{Z}_{M\times N}((i,j))$ and width $d$ (called the annular lattice for short) is defined by

\begin{equation}\label{eqn:6.1}
\mathcal{A}_{M \times N;d}((i,j))=\mathbb{Z}_{(M+2d)\times (N+2d)}((i-d,j-d))\setminus \mathbb{Z}_{M\times N}((i,j)).
\end{equation}
For brevity, let

\begin{equation}\label{eqn:6.2}
\begin{array}{ccc}
\mathcal{A}_{M \times N}((i,j))=\mathcal{A}_{M \times N;2}((i,j)) & \text{and} & \mathcal{A}_{M \times N}=\mathcal{A}_{M \times N}((0,0)).
\end{array}
\end{equation}

The hole-filling condition is defined as follows.

\begin{definition}
\label{definition:6.1}
For $\mathcal{B}\subset\Sigma_{2\times 2}(p)$ and $k\geq 2$, $\mathcal{B}$ satisfies the $k$ hole-filling condition ($($HFC$)_{k}$) with size $(M,N)$, $M,N\geq 2k-3$, if for every $\mathcal{B}$-admissible pattern $U$ on $\mathcal{A}_{M \times N}$ that can be extended to $\mathcal{A}_{(M+4-2k) \times (N+4-2k);k}((k-2,k-2))$ using the local patterns in $\mathcal{B}$, $U$ can completely fill its hole using the patterns in $\mathcal{B}$; see Fig. 5.4. $($HFC$)_{2}$ is also called the hole-filling condition (HFC).

\end{definition}

\begin{equation*}
\begin{array}{ccccc}
\hspace{-1.2cm}\psfrag{a}{{\tiny  $\mathcal{A}_{(M+4-2k) \times (N+4-2k);k}((k-2,k-2))$}}
\psfrag{b}{{\tiny $\mathcal{A}_{M \times N}$}}
\psfrag{k}{{\footnotesize $k$}}
\psfrag{m}{{\footnotesize $M$}}
\psfrag{n}{{\footnotesize $N$}}
\includegraphics[scale=0.7]{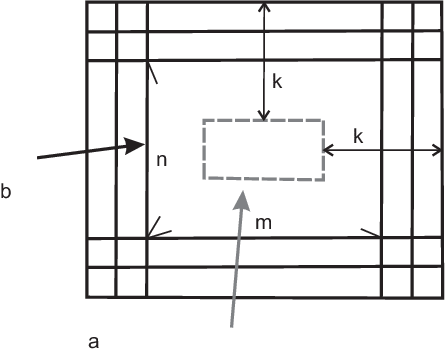} & &  & & \hspace{0.5cm}
\psfrag{a}{{\tiny $$}}
\psfrag{b}{{\tiny $$}}
\psfrag{c}{{\tiny $\xi_{1}$}}
\psfrag{d}{{\tiny $\eta_{1}$}}
\psfrag{e}{{\tiny $\xi_{2}$}}
\psfrag{f}{{\tiny $\eta_{2}$}}
\psfrag{g}{{\tiny $\xi_{N_{1}}$}}
\psfrag{h}{{\tiny $\eta_{N_{1}}$}}
\psfrag{j}{{\tiny  $$}}
\psfrag{k}{{\tiny $$}}
\psfrag{m}{{\tiny $M$}}
\psfrag{n}{{\tiny $N$}}
\psfrag{o}{{\tiny $i_{1}$}}
\psfrag{p}{{\tiny $i_{M}$}}
\psfrag{q}{{\tiny $j_{1}$}}
\psfrag{r}{{\tiny $j_{M}$}}
\psfrag{s}{{\tiny $\cdots$}}
\psfrag{t}{{\tiny $\vdots$}}
\psfrag{u}{{\tiny $(3)$}}
\psfrag{v}{{\tiny $(2)$}}
\psfrag{w}{{\tiny $(4)$}}
\psfrag{x}{{\tiny $(1)$}}
\includegraphics[scale=1.0]{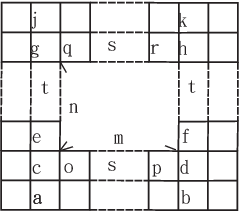}
\\ & & & & \\
\text{Figure 5.4.} & & & & \text{Figure 5.5.}
\end{array}
\end{equation*}

Notably, the $k$ hole-filling condition (HFC$)_{k}$, $k\geq 3$, is weaker than HFC. $($HFC$)_{k}$ can be expressed in terms of the horizontal transition matrices $\mathbb{H}_{n}$ and the connecting operators $S_{m;\alpha,\beta}$ and $W_{m;\alpha,\beta}$. Therefore, the condition $($HFC$)_{k}$ can be easily checked, especially using computer programs. The following theorem concerns only the case in which $\mathcal{B}$ satisfies HFC when $\mathbb{H}_{2}$ and $\mathbb{V}_{2}$ are non-degenerated; for brevity, the general case in which $\mathcal{B}$ satisfies $($HFC$)_{k}$, $k\geq 3$, is omitted.

\begin{theorem}
\label{theorem:6.2}
Given $\mathcal{B}\subset\Sigma_{2\times 2}(p)$, suppose $\mathbb{H}_{2}$ and $\mathbb{V}_{2}$ are non-degenerated. For $M,N\geq 1$, $\mathcal{B}$ satisfies HFC with size $(M,N)$ if and only if for $\alpha_{1},\alpha_{2},\cdots ,\alpha_{N+2}\in\{1,2,\cdots,p^{2}\}$,
%

\begin{equation}\label{eqn:6.5}
S_{M+1;\alpha_{1},\alpha_{2}}S_{M+1;\alpha_{2},\alpha_{3}}\cdots S_{M+1;\alpha_{N+1},\alpha_{N+2}}>0.
\end{equation}
\end{theorem}

\textit{Proof.}
Since $\mathbb{H}_{2}$ and $\mathbb{V}_{2}$ are non-degenerated, from Theorem \ref{theorem:3.5}, $\mathcal{B}$ is rectangle-extendable. Let $N_{1}=N+2$. For any $i_{k}$, $j_{k}$, $\xi_{l}$ and $\eta_{l}\in \mathcal{S}_{p}$, $1\leq k\leq M$ and $1\leq l\leq N_{1}$, a $\mathcal{B}$-admissible pattern can be produced on $\mathcal{A}_{M\times N}$; see Fig. 5.5.
Therefore, by the construction of connecting operators, (\ref{eqn:6.5}) is equivalent to the condition for filling the hole of size $(M,N)$. The proof is complete. \hspace{0.5cm} $\square$
\medbreak

Now, the following theorem provides sufficient conditions for strong specification.

\begin{theorem}
\label{theorem:6.4}
Given $\mathcal{B}\subset\Sigma_{2\times 2}(p)$, if there exists $k\geq 2$ such that
\begin{enumerate}

\item[(i)] $r(\mathbb{H}_{k})=c(\mathbb{H}_{k})$ and $r(\mathbb{V}_{k})=c(\mathbb{V}_{k})$,

\item[(ii)] $\mathcal{B}$ satisfies  $($HFC$)_{k}$ with size $(M,N)$ for some $M,N\geq 2k-3$, and

\item[(iii)] $\mathbb{H}_{k}$ is weakly $(M-2k+5)$-primitive and $\mathbb{V}_{k}$ is weakly $(N-2k+5)$-primitive,
\end{enumerate}
then $\Sigma(\mathcal{B})$ has strong specification.
\end{theorem}

\textit{Proof.}
Let $M'=M-k+4$ and $N'=N-k+4$. First, define the lattice $\mathbb{L}_{g;k}=\mathbb{L}_{g;k}(M,N)$, which is like the grid on a checkerboard with line width $k$ and $(M+4-2k)\times (N+4-2k)$ blank spaces, as

\begin{equation*}
\mathbb{L}_{g;k}=\underset{i,j\in\mathbb{Z}}{\bigcup}\mathcal{A}_{(M+4-2k) \times (N+4-2k);k}((iM'+k-2,jN'+k-2)).
\end{equation*}
Denote the lattice of blank spaces on the checkerboard by

\begin{equation*}
\mathbb{L}_{b;k}=\mathbb{Z}^{2}\setminus \mathbb{L}_{g;k}.
\end{equation*}

\begin{equation*}
\begin{array}{l}
\psfrag{a}{$\cdots$}
\psfrag{b}{$\vdots$}
\psfrag{c}{$(0,0)$}
\includegraphics[scale=0.4]{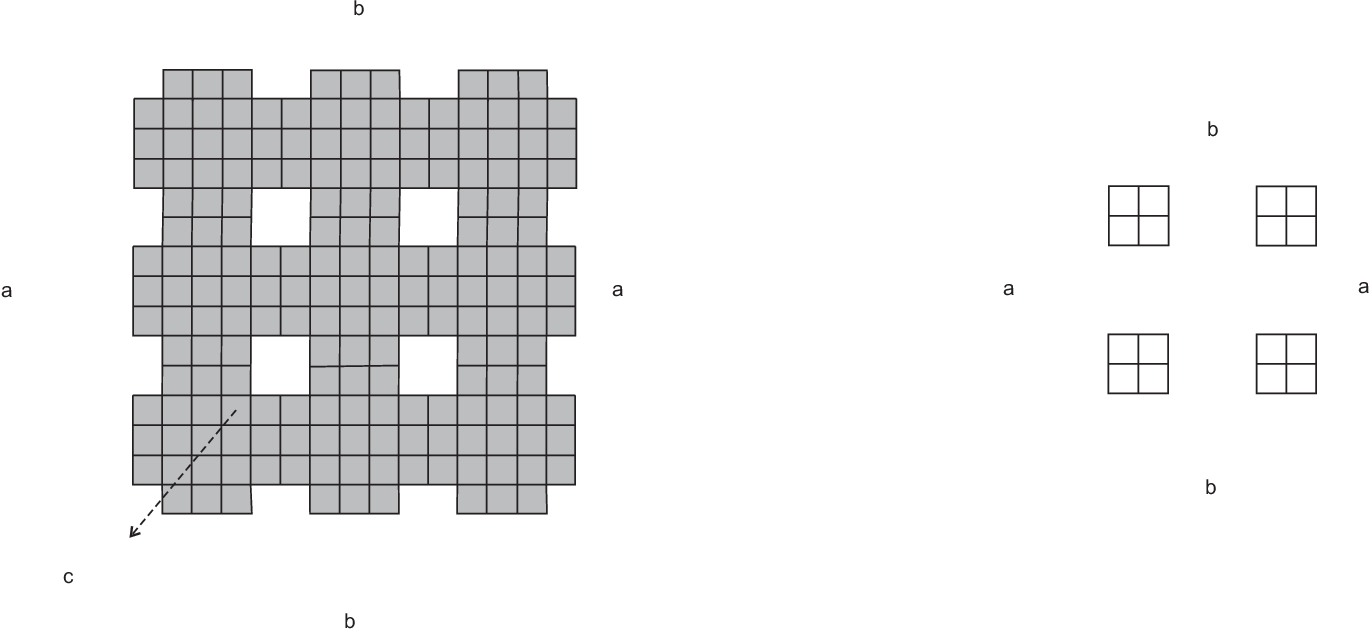}
\\
\\
\text{(a) The shadowed lattice is } \mathbb{L}_{g;3}. \hspace{1.2cm}
\text{(b) } \text{The white lattice is }\mathbb{L}_{b;3}.\\
\ \hspace{0.5cm}\text{for }M=N=4. \hspace{3.3cm}\text{ for }M=N=4.
\end{array}
\end{equation*}
\begin{equation*}
\text{Figure 5.6.}
\end{equation*}

For $i,j\in\mathbb{Z}$, define

\begin{equation*}
\left\{
\begin{array}{l}
\mathbb{L}(i,j)=\mathbb{L}_{k;M,N}(i,j)=\mathbb{Z}_{M'\times N'}\left(\left(iM'-2,jN'-2\right)\right) \\
\\
\widehat{\mathbb{L}}(i,j)=\widehat{\mathbb{L}}_{k;M,N}(i,j)=\mathbb{Z}_{(M+4)\times(N+4)}\left(\left(iM'-2,jN'-2\right)\right).
\end{array}
\right.
\end{equation*}

\begin{equation*}
\begin{array}{c}
\psfrag{a}{{\footnotesize $\left(5i-2,5j-2\right)$}}
\psfrag{b}{ $\mathbb{L}_{3;4,4}(i,j)=$}
\psfrag{c}{ $\widehat{\mathbb{L}}_{3;4,4}(i,j)=$}
\psfrag{d}{and}
 \includegraphics[scale=0.6]{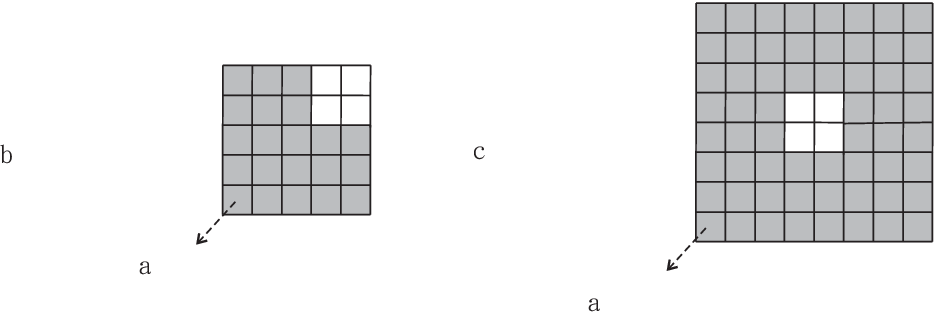}
\end{array}
\end{equation*}

\begin{equation*}
\text{Figure 5.7. The lattices }\mathbb{L}_{3;4,4}(i,j)\text{ and }\widehat{\mathbb{L}}_{3;4,4}(i,j).
\end{equation*}

Clearly, $\widehat{\mathbb{L}}(i,j)\supset\mathbb{L}(i,j)$ and $\mathbb{L}(i_{1},j_{1})\bigcap\mathbb{L}(i_{2},j_{2})=\emptyset$ if $(i_{1},j_{1})\neq (i_{2},j_{2})$. Then, $\mathbb{Z}^{2}$ lattice can be decomposed into disjoint sublattices:

\begin{equation*}
\mathbb{Z}^{2}  =\underset{i,j\in\mathbb{Z}}{\bigcup}\mathbb{L}(i,j).
\end{equation*}

Take
\begin{equation}\label{eqn:6.7}
\bar{d}=3\sqrt{(M')^{2}+(N')^{2}}.
\end{equation}
Let $R_{1},R_{2}\subset \mathbb{Z}^{2}$ with $d(R_{1},R_{2})\geq\bar{d}$.
For any $U_{l}=\Pi_{R_{l}}(W_{l})$ with $W_{l}\in\Sigma(\mathcal{B})$, $l=1,2$, let

\begin{equation*}
R_{l}'=\underset{(i',j')\in\mathbb{Z}_{2\times 2}((i-1,j-1))}{\underset{R_{l}\bigcap \mathbb{L}(i,j)\neq\emptyset}{\bigcup}}\widehat{\mathbb{L}}(i',j')
\end{equation*}
for $l=1,2$. Hence, $U_{l}$ can be extended as $U_{l}'=\Pi_{R_{l}'}(W_{l})$, $l=1,2$.
Clearly, for $l=1,2$,

\begin{equation}\label{eqn:6.7-1}
\text{if } (i,j)\in R_{l}, \text{ then } \mathbb{Z}_{(2k+1)\times(2k+1)}((i-k,j-k))\subseteq R_{l}'.
\end{equation}

From (\ref{eqn:6.7}), it can be verified that $\widehat{\mathbb{L}}(i,j)\bigcap R_{1}'\neq\emptyset$ and $\widehat{\mathbb{L}}(i,j)\bigcap R_{2}'\neq\emptyset$ never both occur for all $(i,j)\in\mathbb{Z}^{2}$.

Now, from conditions (i) and (iii), there exists a $\mathcal{B}$-admissible pattern $U''$ on $R_{1}'\bigcup R_{2}' \bigcup \mathbb{L}_{g;k}$ such that $U''\mid_{R_{l}'}=U_{l}'$, $i=1,2$.
Clearly, $\mathbb{Z}^{2}\setminus\left(R_{1}'\bigcup R_{2}' \bigcup \mathbb{L}_{c}\right)$ is the union of the discrete $(M+4-2k)\times (N+4-2k)$ rectangular lattices.

Hence, from (\ref{eqn:6.7-1}) and condition (2), there exists $W\in\Sigma(\mathcal{B})$ such that $W\mid_{R_{i}}=U_{i}$ for $i=1,2$. Notably, in general, $W\mid_{R_{1}'\bigcup R_{2}' \bigcup \mathbb{L}_{c} }$ is not equal to $U''$ since condition (2) may change the colors on the boundary of $R_{1}'\bigcup R_{2}' \bigcup \mathbb{L}_{c}$ with width $k-2$. Therefore, $\Sigma(\mathcal{B})$ has strong specification. The proof is complete. \hspace{0.5cm} $\square$
\medbreak

Theorem \ref{theorem:6.4} clearly applies to certain non-degenerated cases, including the Golden Mean shift. The Golden Mean shift is known to have safe symbol $0$ and strong specification; see \cite{50}. The Golden Mean shift can also be shown to satisfy HFC with size $(1,1)$ and to have strong specification by Theorem \ref{theorem:6.4}.

The following well-known example of Burton and Steif \cite{12,13} is introduced for the further application of Theorem \ref{theorem:6.4}. This example is closely related to the ferromagnetic Ising model in statistical physics.

\begin{example}
\label{example:6.6}
Consider the color set $\mathcal{S}_{4}'=\{-2,-1,1,2\}$. The rule of $\mathbf{X}_{BS}\subseteq \mathcal{S}_{4}'^{\mathbb{Z}^{2}}$ is that a negative is disallowed to sit to a positive unless they are both $\pm 1$. To fit  $\mathcal{S}_{4}'$ to the color set $\mathcal{S}_{4}=\{0,1,2,3\}$ used in this work, $-2$, $-1$, $1$ and $2$ are replaced with $0$, $1$, $2$ and $3$, respectively.
That $\mathcal{B}_{BS}$ satisfies HFC with size $(2,2)$ can be proven and the details are omitted. Therefore, by Theorem \ref{theorem:6.4}, $\Sigma(\mathcal{B}_{BS})$ has strong specification.

\end{example}

For $p=2$, the size $(M,N)$ of the hole-filling condition can be larger, as in the following example.

\begin{example}
\label{example:6.13}
From Theorem \ref{theorem:6.2}, it can be verified that
\begin{equation*}
\mathbb{H}_{2}(\mathcal{B}_{1})=
\left[
\begin{array}{cccc}
1 & 1 & 1 & 1 \\
1 & 1 & 0 & 1 \\
1& 0 & 1 & 1 \\
1 & 1 & 1 & 0
\end{array}
\right]
\end{equation*}
satisfies HFC with size $(3,3)$ and

\begin{equation*}
\mathbb{H}_{2}(\mathcal{B}_{2})=
\left[
\begin{array}{cccc}
1 & 1 & 1 & 1 \\
1 & 0 & 1 & 1 \\
1& 1 & 1 & 1 \\
0 & 1 & 1 & 1
\end{array}
\right]
\end{equation*}
satisfies HFC with size $(4,4)$. Therefore, by Theorem \ref{theorem:6.4}, both can be shown to have strong specification. Notably, they do not have a safe symbol $0$ or $1$. The details are omitted.

\end{example}

The following example concerns the Diagonally Restricted Golden Mean, which does not satisfy HFC but does satisfy $($HFC$)_{3}$.

\begin{example}
\label{example:6.15}(Diagonally Restricted Golden Mean) Consider $\mathcal{S}_{2}=\{0,1\}$ and

\begin{equation*}
\mathbb{H}_{2}(\mathcal{B}_{s})=\mathbb{V}_{2}(\mathcal{B}_{s})=
\left[
\begin{array}{cccc}
1 & 1 & 1 & 0 \\
1 & 0 & 0 & 0 \\
1& 0 & 0 & 0 \\
0 & 0 & 0 & 0
\end{array}
\right].
\end{equation*}

Since \psfrag{c}{{\scriptsize $1$}}
\psfrag{b}{{\scriptsize $0$}}
\psfrag{d}{{\scriptsize $1$}}
\psfrag{e}{}
\includegraphics[scale=0.5]{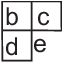} (or \psfrag{c}{{\scriptsize $1$}}
\psfrag{e}{{\scriptsize $0$}}
\psfrag{d}{{\scriptsize $1$}}
\psfrag{b}{}
\includegraphics[scale=0.5]{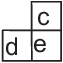}) is forbidden, $\mathcal{B}_{s}$ is easily seen not to satisfy HFC.

The fact that $\mathcal{B}_{s}$ satisfies $($HFC$)_{3}$ with size $(3,3)$ can be verified. Clearly,
$r(\mathbb{H}_{3})=c(\mathbb{H}_{3})$, and $\mathbb{H}_{3}=\mathbb{V}_{3}$ is $2$-primitive. Therefore, by  Theorem \ref{theorem:6.4}, $\Sigma (\mathcal{B}_{s})$ has strong specification. Notably, $0$ is known to be a safe symbol and strong specification follows immediately.

\end{example}

\newpage

\begin{equation*}
\text{\textbf{Appendix}}
\end{equation*}

This appendix recalls the various mixing properties that were described in the Introduction; see Boyle et al. \cite{11}.

\textbf{Definition A.1.}\emph{ Suppose $\Sigma$ is a $\mathbb{Z}^{2}$ shift.}
\begin{itemize}
\item[(i)]\emph{$\Sigma$ has the uniform filling property (UFP) if a number $M(\Sigma)\geq 1$ exists such that for any two allowable patterns $U_{1}\in\Pi_{R_{1}}(\Sigma)$ and $U_{2}\in\Pi_{R_{2}}(\Sigma)$ with $d(R_{1},R_{2})\geq M$, where $R_{1}=\mathbb{Z}_{m\times n}((i,j))$, $m,n\geq 1$ and $(i,j)\in\mathbb{Z}^{2}$, and  $R_{2}\subset\mathbb{Z}^{2}$, there exists a global pattern $W\in\Sigma$ with $\Pi_{R_{1}}(W)=U_{1}$ and $\Pi_{R_{2}}(W)=U_{2}$.}

\item[(ii)]\emph{$\Sigma$ is strongly irreducible if a number $M(\Sigma)\geq 1$ exists such that for any two allowable patterns $U_{1}\in\Pi_{R_{1}}(\Sigma)$ and $U_{2}\in\Pi_{R_{2}}(\Sigma)$ with $d(R_{1},R_{2})\geq M$, where $R_{1}\subset\mathbb{Z}^{2}$ is finite and  $R_{2}\subset\mathbb{Z}^{2}$, there exists a global pattern $W\in\Sigma$ with $\Pi_{R_{1}}(W)=U_{1}$ and $\Pi_{R_{2}}(W)=U_{2}$.}

\item[(iii)]\emph{$\Sigma$ is corner gluing if a number $M(\Sigma)\geq 1$ exists such that for any two allowable patterns $U_{1}\in\Pi_{R_{1}}(\Sigma)$ and $U_{2}\in\Pi_{R_{2}}(\Sigma)$ with $d(R_{1},R_{2})\geq M$, where $R_{1}=\mathbb{Z}_{m\times n}((i,j))$, $m,n\geq 1$ and $(i,j)\in\mathbb{Z}^{2}$, and  $R_{2}=\mathbb{Z}_{m_{1}\times n_{1}}((i+m-m_{1},j+n-n_{1}))\setminus \mathbb{Z}_{m_{2}\times n_{2}}((i+m-m_{2},j+n-n_{2}))$, $m_{1}>m_{2}\geq m+M$ and $n_{1}>n_{2}\geq n+M$, there exists a global pattern $W\in\Sigma$ with $\Pi_{R_{1}}(W)=U_{1}$ and $\Pi_{R_{2}}(W)=U_{2}$.}

\item[(iv)]\emph{$\Sigma$ is block gluing if a number $M(\Sigma)\geq 1$ exists such that for any two allowable patterns $U_{1}\in\Pi_{R_{1}}(\Sigma)$ and $U_{2}\in\Pi_{R_{2}}(\Sigma)$ with $d(R_{1},R_{2})\geq M$, where $R_{1}=\mathbb{Z}_{m_{1}\times n_{1}}((i_{1},j_{1}))$ and  $R_{2}=\mathbb{Z}_{m_{2}\times n_{2}}((i_{2},j_{2}))$, $m_{l},n_{l}\geq1 $ and $(i_{l},j_{l})\in\mathbb{Z}^{2}$, $l\in\{1,2\}$, there exists a global pattern $W\in\Sigma$ with $\Pi_{R_{1}}(W)=U_{1}$ and $\Pi_{R_{2}}(W)=U_{2}$.}

\end{itemize}
Notably, (i)$\sim$(iv) were introduced in \cite{11,29-1,47,50}.

Their significance in the classification of mixing properties is discussed as follows. Boyle et al. \cite{11} discussed various mixing properties, including strong irreducibility, the uniform filling property (UFP), corner gluing, block gluing and topological mixing.
Figure A.1 presents the range of (HFC$)_{k}$ and these mixing properties for $\mathbb{Z}^{2}$ shifts of finite type. For brevity, the following notation is used.
\begin{enumerate}
\item[(a)]: $\mathcal{B}$ satisfies (HFC$)_{k}$ and conditions (i) and (iii) of Theorem \ref{theorem:1.2},

\item[(b)]: $\Sigma(\mathcal{B})$ has strong specification,


\item[(c)]:  $\Sigma(\mathcal{B})$ has the UFP,

\item[(d)]:  $\Sigma(\mathcal{B})$ is corner gluing,

\item[(e)]:  $\Sigma(\mathcal{B})$ is block gluing,

\item[(f)]:  $\Sigma(\mathcal{B})$ is topologically mixing.
\end{enumerate}

\begin{equation*}
\psfrag{a}{{\footnotesize (a)}}
\psfrag{b}{{\footnotesize (b)}}
\psfrag{c}{{\footnotesize (c)}}
\psfrag{d}{{\footnotesize (d)}}
\psfrag{e}{{\footnotesize (e)}}
\psfrag{f}{{\footnotesize (f)}}
\includegraphics[scale=0.4]{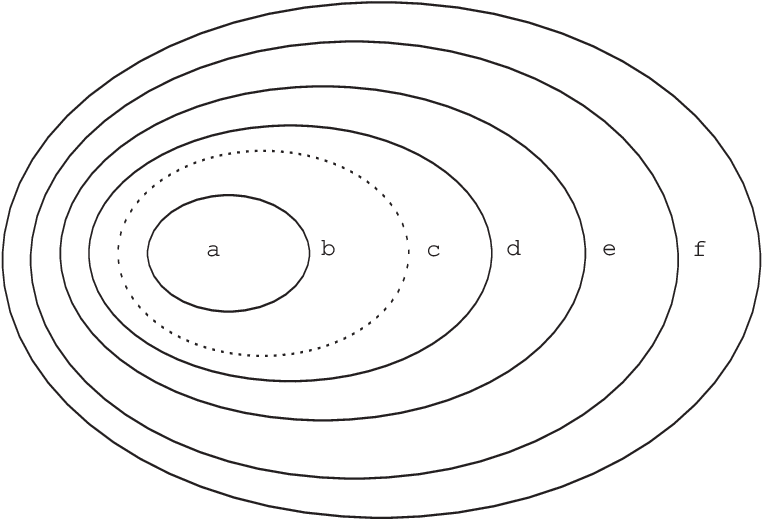}
\end{equation*}
\begin{equation*}
\text{Figure A.1.}
\end{equation*}

Notably, the solid line indicates that the inward property is strictly stronger than the outward one; the dotted line indicates that the inward property is stronger than or equivalent to the outward one. The examples for $(d)\nRightarrow (c)$ and $(e)\nRightarrow (d)$ were given by Boyle et al. \cite{11}. (b) and (c) are not the same for general subshifts \cite{31-1}; the equivalence is still open for shifts of finite type.

With reference to Fig. A.1, Theorem \ref{theorem:1.1} and the primitive results in Section 4 ensure that the weakest case--topological mixing--holds. Theorem \ref{theorem:1.2} ensures that the strongest case--strong specification--holds.

\end{document}